\documentclass[a4paper,10pt,BCOR1cm,parskip]{scrartcl}

\usepackage{german}
\usepackage[iso]{umlaute}
\usepackage{theorem}
\usepackage{amsmath}
\usepackage{amssymb}
\usepackage{makeidx}
\usepackage[matrix,arrow]{xy}
\usepackage{epsfig}
\usepackage{boxedminipage}
\usepackage{epic}
\usepackage{longtable}
\usepackage{ifthen}
\usepackage{tabularx}
\usepackage{curves}
\usepackage{parskip}

\newcommand{\R}{\mathrm{I\!R}}
\newcommand{\N}{\mathrm{I\!N}}
\newcommand{\HH}{\mathrm{I\!H}}

\newcommand{\Z}{\mathchoice {\hbox{$\sf\textstyle Z\kern-0.4em
Z$}}{\hbox{$\sf\textstyle Z\kern-0.4em Z$}}{\hbox{$\sf\scriptstyle
Z\kern-0.3em Z$}}{\hbox{$\sf\scriptscriptstyle Z\kern-0.2em Z$}}}

\newcommand{\Q}{\mathchoice {\setbox0=\hbox{$\displaystyle\rm
Q$}\hbox{\raise0.15\ht0\hbox to0pt{\kern0.4\wd0\vrule
height0.8\ht0\hss}\box0}}{\setbox0=\hbox{$\textstyle\rm
Q$}\hbox{\raise0.15\ht0\hbox to0pt{\kern0.4\wd0\vrule
height0.8\ht0\hss}\box0}}{\setbox0=\hbox{$\scriptstyle\rm
Q$}\hbox{\raise0.15\ht0\hbox to0pt{\kern0.4\wd0\vrule
height0.7\ht0\hss}\box0}}{\setbox0=\hbox{$\scriptscriptstyle\rm
Q$}\hbox{\raise0.15\ht0\hbox to0pt{\kern0.4\wd0\vrule
height0.7\ht0\hss}\box0}}}

\newcommand{\id}{\mathrm{id}}
\newcommand{\GL}{\mathrm{GL}}

\newcommand{\SO}{\mathrm{SO}}
\newcommand{\End}{\mathrm{End}}

\newcommand{\rk}{\mathrm{rk}}

\newcommand{\eps}{\varepsilon}
\newcommand{\vi}{\varphi}
\newcommand{\vkap}{\varkappa}

\newcommand{\qmq}[1]{\quad\mbox{#1}\quad}
\newcommand{\Menge}[2]{\{\,#1\,|\,#2\,\}}

\renewcommand{\bigoplus}{\mathop{\bigcirc
  \raisebox{-0.22em}{\hskip-0.53em\hbox{\vrule height2.08ex width0.04em}
  \raisebox{ 0.48em}{\hskip-0.75em\hbox{\vrule height0.04em width 0.8em}}
  \hskip- 0.2em}}}

\newcommand {\g}[2]{\langle #1,#2\rangle}

\newcommand {\operp}{\mathbin{\mbox{$\ominus\raisebox{2.9pt}
 {\hskip-0.42em\hbox{\vrule height0.7ex width0.02em}\hskip0.42em }$}}}

\newcommand{\Ug}{\mathrm{U}}

\newcommand{\RE}{\mathop{\mathrm{Re}}\nolimits}
\newcommand{\IM}{\mathop{\mathrm{Im}}\nolimits}

\newcommand{\Eig}{\mathop{\mathrm{Eig}}\nolimits}
\newcommand{\Aut}{\mathop{\mathrm{Aut}}\nolimits}

\newcommand{\ad}{\mathop{\mathrm{ad}}\nolimits}
\newcommand{\Ad}{\mathop{\mathrm{Ad}}\nolimits}
\newcommand{\Exp}{\mathop{\mathrm{Exp}}\nolimits}
\newcommand{\Fix}{\mathop{\mathrm{Fix}}\nolimits}

\newcommand{\tr}{\mathop{\mathrm{tr}}\nolimits}

\newcommand{\Geins}{\mathrm{G1}}
\newcommand{\Gzwei}{\mathrm{G2}}
\newcommand{\Gdrei}{\mathrm{G3}}
\newcommand{\Peins}{\mathrm{P1}}
\newcommand{\Pzwei}{\mathrm{P2}}
\newcommand{\Atyp}{\mathrm{A}}
\newcommand{\Ieins}{\mathrm{I1}}
\newcommand{\Izwei}{\mathrm{I2}}

\newcommand{\A}{\mathfrak{A}}
\newcommand{\liea}{\mathfrak{a}}

\newcommand{\lieg}{\mathfrak{g}}

\newcommand{\liek}{\mathfrak{k}}

\newcommand{\liem}{\mathfrak{m}}
\newcommand{\lieo}{\mathfrak{o}}

\newcommand{\RP}{\ensuremath{\R\mathrm{P}}}
\newcommand{\CP}{\ensuremath{\C\mathrm{P}}}

\newcommand{\CQ}{\ensuremath{\C\mathrm{Q}}}

\newcommand{\Geo}{\mathrm{Geo}}

\newcommand{\scrH}{\mathcal{H}}
\newcommand{\scrI}{\mathcal{I}}
\newcommand{\scrR}{\mathcal{R}}
\newcommand{\bbS}{\mathbb{S}}
\newcommand{\bbT}{\mathbb{T}}
\newcommand{\bbV}{\mathbb{V}}

\newcommand{\Sph}{\bbS}

\newcommand{\gC}[2]{\langle #1,#2\rangle_{\C}}

\newcommand{\beweis}{\begingroup\footnotesize \emph{Proof. }}
\newcommand{\beweisende}{\strut\hfill $\Box$\par\medskip\endgroup}

\newcommand{\nsp}[2]{\perp_{#2}\nobreak\!\nobreak\!\nobreak #1}		% "normal space in p"
\newcommand{\unsp}[2]{\perp^1_{#2}\nobreak\!\nobreak\!\nobreak #1}	% "unit normal space in p"

\newcommand{\wt}{\widetilde}
\newcommand{\wh}{\widehat}

\newcommand{\C}{\mathchoice {\setbox0=\hbox{$\displaystyle\rm
C$}\hbox{\hbox to0pt{\kern0.4\wd0\vrule
height0.95\ht0\hss}\box0}}{\setbox0=\hbox{$\textstyle\rm C$}\hbox{\hbox
to0pt{\kern0.4\wd0\vrule
height0.95\ht0\hss}\box0}}{\setbox0=\hbox{$\scriptstyle\rm C$}\hbox{\hbox
to0pt{\kern0.4\wd0\vrule
height0.95\ht0\hss}\box0}}{\setbox0=\hbox{$\scriptscriptstyle\rm
C$}\hbox{\hbox to0pt{\kern0.4\wd0\vrule height0.95\ht0\hss}\box0}}}

\theoremstyle{plain} % margin, change

\theorembodyfont{\rmfamily\mdseries\itshape}
\newtheorem{Def}{Definition}[section]
\newtheorem{Prop}[Def]{Proposition}
\newtheorem{Theorem}[Def]{Theorem}
\newtheorem{Lemma}[Def]{Lemma}

\theorembodyfont{\rmfamily\mdseries\upshape}

\newtheorem{Examples}[Def]{Examples}
\newtheorem{Remark}[Def]{Remark}

\begin{document}
\selectlanguage{english}

\title{Totally geodesic submanifolds of the complex quadric}
\author{Sebastian Klein}
\date{March 7, 2006}
\maketitle

\abstract{\textbf{Abstract.} In this article, relations between the root space decomposition of a Riemannian symmetric space of compact type and the
root space decompositions of its totally geodesic submanifolds (symmetric subspaces) are described. These relations provide an
approach to the classification of totally geodesic submanifolds in Riemannian symmetric spaces. 
In this way a classification of the totally geodesic submanifolds in the complex quadric \,$Q^m := \SO(m+2)/(\SO(2) \times \SO(m))$\,
is obtained. It turns out that the earlier classification of totally geodesic submanifolds of \,$Q^m$\, 
by \textsc{Chen} and \textsc{Nagano} is incomplete: in particular a type of submanifolds which are isometric to \,$2$-spheres of radius \,$\tfrac{1}{2}\sqrt{10}$\,, and 
which are neither complex nor totally real in \,$Q^m$\,, is missing.}

\section{Introduction}
\label{Se:intro}

This article is concerned with the study of totally geodesic submanifolds in Riemannian symmetric spaces of compact type. Its objective is two-fold:
First, we describe general relations between the roots and root spaces of such a symmetric space, and the roots resp.~root spaces of its totally geodesic
submanifolds. Second, we apply these results to obtain a classification of the totally geodesic submanifolds in the complex quadric
\,$Q^m = \SO(m+2)/(\SO(2)\times \SO(m))$\,. 

\bigskip

It should be mentioned that already \textsc{Chen} and \textsc{Nagano} gave a classification of the totally geodesic submanifolds of the complex quadric
by ``ad-hoc methods'' in \cite{Chen/Nagano:totges1-1977}. However that paper contains (besides some inaccuracies which are easily resolved)
a more serious mistake, which causes two types of totally geodesic submanifolds to be missed. The submanifolds of the first of these two types are
isometric to \,$\CP^1 \times \RP^1$\,; as it is explained in Section~\ref{Se:tgsub}, 
their existence can be derived from the fact that \,$Q^2$\, is holomorphically isometric to \,$\CP^1 \times \CP^1$\, (via the Segre embedding).
The manifolds of the second type are isometric to $2$-spheres of radius \,$\tfrac{1}{2}\sqrt{10}$\,. They are neither complex nor
totally real submanifolds of \,$Q^m$\,, and they are remarkable insofar as their geodesic diameter \,$\tfrac{\pi}{2}\sqrt{10}$\, is strictly
larger than the geodesic diameter \,$\tfrac{\pi}{\sqrt{2}}$\, of \,$Q^m$\,.

In \cite{Chen/Nagano:totges2-1978}, Chen and Nagano introduced their \,$(M_+,M_-)$-method for the classification of totally geodesic submanifolds
in Riemannian symmetric spaces of compact type, and via this method they again give a classification of the totally geodesic submanifolds in
rank 2 symmetric spaces. However, the totally geodesic submanifolds of \,$Q^m$\, which were missing from \cite{Chen/Nagano:totges1-1977} are
also missing here.

\bigskip

The approach to the classification of totally geodesic submanifolds taken here is as follows:
It is well-known that in any symmetric space \,$M=G/K$\,, for given \,$p \in M$\, and \,$U \subset T_pM$\,, there exists a totally geodesic
submanifold \,$M'$\, of \,$M$\, with \,$p \in M'$\, and \,$T_pM' = U$\, if and only if \,$U$\, is curvature-invariant (i.e.~if \,$R(u,v)w \in U$\,
holds for every \,$u,v,w \in U$\,, where \,$R$\, is the curvature tensor of \,$M$\,). Moreover, if we consider the decomposition \,$\lieg = \liek \oplus \liem$\,
of the Lie algebra \,$\lieg$\, of \,$G$\, induced by the symmetric space structure of \,$M$\, (we then have the canonical isomorphism \,$\tau:\liem \to T_pM$\,), 
then \,$U \subset T_pM$\, is curvature-invariant
if and only if \,$\liem' := \tau^{-1}(U) \subset \liem$\, is a Lie triple system (i.e.~if \,$[[\liem',\liem'],\liem'] \subset \liem'$\,
holds).
Thus, the task of classifying the totally geodesic submanifolds of \,$M$\, splits into two steps: (1) To classify the Lie triple systems
in \,$\liem$\,, and (2) for each of the Lie triple systems \,$\liem'$\, found in the first step, to construct a totally geodesic, connected, complete
submanifold \,$M'$\, of \,$M$\, so that \,$p \in M'$\, and \,$\tau^{-1}(T_pM') = \liem'$\, holds.

In Section~\ref{Se:roots}, Lie triple systems in Riemannian symmetric spaces of compact type are studied with regard to the theory of roots and
root spaces. In particular, relations between the roots and root spaces of a Lie triple system and the roots resp.~root spaces of the 
ambient symmetric space are derived. These relations turn out to be useful for the classification of the Lie triple systems,
at least for the complex quadric, as will be seen.

The remainder of the article is concerned with the application of these general results to the complex quadric \,$Q^m = \SO(m+2)/(\SO(2)\times \SO(m))$\,
(a rank \,$2$\, Riemannian symmetric space); thereby a classification of the totally geodesic submanifolds of \,$Q^m$\, is obtained. In
Section~\ref{Se:Q} we describe some facts regarding the geometry of the complex quadric which are needed for the classification. These
facts are mostly taken from the paper \cite{Reckziegel:quadrik-1995} by \textsc{H.~Reckziegel}; especially the concept of a \CQ-structure
(see Definition~\ref{D:Q:CQ}), which is very useful for the formulation of the classification, was introduced there. 

In Section~\ref{Se:cla} the main result of the present article, the classification of the Lie triple systems of the complex quadric, is obtained,
see Theorem~\ref{T:cla:cla}. The proof of this theorem is based on the combination of the general results on roots and root spaces
of Lie triple systems from Section~\ref{Se:roots} with the specific description of the geometry of the complex quadric given in Section~\ref{Se:Q}.

Finally, in Section~\ref{Se:tgsub} the totally geodesic submanifolds which correspond to the various Lie triple systems described in Theorem~\ref{T:cla:cla}
are described; thereby the classification of the totally geodesic submanifolds of \,$Q^m$\, is completed (see the table given in that section).

\bigskip

The results presented in the present paper were obtained by me in my dissertation under the advisorship of Professor \textsc{H.~Reckziegel}.
I wish to express my sincerest gratitude for his enduring and intensive support. 

\section{The root space decomposition corresponding to a Lie triple system}
\label{Se:roots}

In this section we suppose that \,$M=G/K$\, is any Riemannian symmetric space of compact type. 
We consider the decomposition \,$\lieg = \liek \oplus \liem$\, of the Lie algebra \,$\lieg$\, of \,$G$\, induced by the symmetric
structure of \,$M$\,. Because \,$M$\, is of compact type, the Killing form \,$\vkap: \lieg \times \lieg \to \R,\;(X,Y) \mapsto \tr(\ad(X) \circ \ad(Y))$\,
is negative definite, and therefore \,$\g{\cdot}{\cdot} := -c \cdot \vkap$\, gives rise to a Riemannian metric on \,$M$\, for any \,$c \in \R_+$\,.%
\footnote{The choice of the factor \,$c$\, does not have any geometric significance. The only reason for considering such a factor is to
accommodate the natural Riemannian metric on the complex quadric, see Equation~\eqref{eq:Q:skp-killing} below.}
In the sequel we suppose that \,$M$\, is equipped with such a Riemannian metric.

\begin{Def}
A linear subspace \,$\liem' \subset \liem$\, is called a \emph{Lie triple system} if \,$[\,[\liem',\liem'] \,,\, \liem'\,] \subset \liem'$\, holds.
\end{Def}

As was explained in the introduction, the first step in classifying the totally geodesic submanifolds of \,$M$\,
is to classify the Lie triple systems of \,$\liem$\,. In this section we describe general results concerning the relationship between the roots resp.~root spaces
of \,$M$\, and those of its totally geodesic submanifolds; in the case \,$M=Q^m$\, these results will permit to classify the Lie triple systems of \,$Q^m$\,. 

First we fix notations concerning flat subspaces, the roots and root spaces of \,$M$\, (for the corresponding theory, see for example \cite{Loos:1969-2}, Section~V.2):
A linear subspace \,$\liea \subset \liem$\, is called flat if \,$[\liea,\liea] = \{0\}$\, holds. The maximal dimension of a flat subspace of \,$\liem$\,
is called the rank of \,$M$\, (or of \,$\liem$\,) and is denoted by \,$\rk(M)$\,. The flat subspaces \,$\liea \subset \liem$\, with \,$\dim(\liea) = \rk(M)$\,
are called the Cartan subalgebras (or maximal flat subspaces) of \,$\liem$\,. We now fix a Cartan subalgebra 
\,$\liea \subset \liem$\,. Then we put for any linear form \,$\lambda \in \liea^*$\,
$$ \liem_\lambda := \Menge{\,X \in \liem\,}{\,\forall Z \in \liea: \ad(Z)^2 X = -\lambda(Z)^2X\,} $$
and consider the root system
$$ \Delta(\liem,\liea) := \Menge{\,\lambda \in \liea^*\setminus \{0\}\,}{\,\liem_\lambda \neq \{0\}\,} $$
of \,$\liem$\, with respect to \,$\liea$\,. (The elements of \,$\Delta(\liem,\liea)$\, are called roots of \,$\liem$\, with respect to \,$\liea$\,,
and for \,$\lambda \in \Delta(\liem,\liea)$\,, \,$\liem_\lambda$\, is called the root space corresponding to \,$\lambda$\,.) As is well-known, we have 
\begin{equation}
\label{eq:roots:m0a}
\liem_0 = \liea
\end{equation}
and
\begin{equation}
\label{eq:roots:mdecomp}
\liem = \liea \;\oplus\; \bigoplus_{\lambda \in \Delta_+} \liem_\lambda \;,
\end{equation}
where \,$\Delta_+ \subset \Delta(\liem,\liea)$\, is any system of positive roots, i.e.~we have \,$\Delta_+ \cup (-\Delta_+) = \Delta(\liem,\liea)$\, 
and \,$\Delta_+ \cap (-\Delta_+) = \varnothing$\,.

Let us fix a Lie triple system \,$\liem' \subset \liem$\,. It should be noted that \,$\liem'$\, does not need to be of compact type, and
therefore the usual root theory for symmetric spaces is not applicable to \,$\liem'$\, directly. However, we will now see how the fact that
\,$\liem'$\, is contained in the space \,$\liem$\, of compact type can be used to construct a root space decomposition for \,$\liem'$\,.
We base this construction on the fact (see \cite{Loos:1969-1}, Theorem~IV.1.6, p.~145 and its proof) that \,$\liem'$\, can be decomposed into Lie triple systems 
$$ \liem' = \liem'_{fl} \oplus \liem'_c \oplus \liem'_{nc} \;, $$
where \,$\liem'_{fl}$\, is flat, \,$\liem'_{c}$\, corresponds to a symmetric space of compact type and \,$\liem'_{nc}$\, corresponds to a symmetric
space of non-compact type; moreover we have
\begin{equation}
\label{eq:roots:km'0}
[[\liem',\liem'],\liem'_{fl}] = \{0\}
\end{equation}
and
\begin{equation}
\label{eq:roots:mm'0}
[\liem',\liem'_{fl}] = \{0\}\; . 
\end{equation}

In fact, in the present situation \,$\liem'_{nc} = \{0\}$\, holds, as the following argument shows: Because the Riemannian symmetric space \,$M$\, is of compact type,
its sectional curvature is \,$\geq 0$\,. The totally geodesic submanifold \,$M'$\, of \,$M$\, corresponding to \,$\liem'$\, 
therefore also has sectional curvature \,$\geq 0$\,. This means in particular that the Ricci curvature form of \,$M'$\,
is positive semi-definite. On the other hand, if \,$\liem'_{nc}$\, were non-zero, it would correspond to a symmetric subspace of \,$M'$\, of non-compact type,
whose Ricci curvature form would be negative definite, in contradiction to the preceding statement. Thus we in fact have
\begin{equation}
\label{eq:roots:m'flc}
\liem' = \liem'_{fl} \oplus \liem'_c \; .
\end{equation}

Note that because \,$\liem'_c$\, is of compact type, the usual concepts of rank, Cartan subalgebras, roots and root systems are applicable to \,$\liem'_c$\,.

In analogous application of the usual concepts we call the maximal dimension of a flat subspace of \,$\liem'$\, the rank of \,$\liem'$\,, 
which we denote by \,$\rk(\liem')$\,.
Moreover, we call a flat subspace \,$\liea' \subset \liem'$\, a Cartan subalgebra of \,$\liem'$\, if \,$\dim(\liea') = \rk(\liem')$\, holds.

%\begin{Prop}
%\label{P:roots:flat}
%A linear subspace \,$\liea' \subset \liem'$\, is a Cartan subalgebra of \,$\liem'$\, if and only if \,$\liea' = \liem'_{fl} \oplus \liea'_c$\, holds
%with some Cartan subalgebra \,$\liea'_c$\, of \,$\liem'_c$\,. We have \,$\rk(\liem') = \dim(\liem_{fl}') + \rk(\liem_{c}')$\,. 
%\end{Prop}
%
%\beweis
%We first note that because of Equation~\eqref{eq:roots:mm'0} we have in any case
%\begin{equation}
%\label{eq:roots:flat:a'liem'fl}
%[\liea',\liem'_{fl}] = \{0\} \; . 
%\end{equation}
%
%Now suppose that \,$\liea'$\, is a Cartan subalgebra of \,$\liem'$\,. 
%Because \,$\liea'$\, and \,$\liem'_{fl}$\, are flat, Equation~\eqref{eq:roots:flat:a'liem'fl} shows that also \,$\liea' + \liem'_{fl}$\, is flat,
%and therefore it follows by the maximality of \,$\liea'$\, that \,$\liem'_{fl} \subset \liea'$\, holds. The decomposition~\eqref{eq:roots:m'flc} then shows that
%there exists a linear subspace \,$\liea'_c \subset \liem'_c$\, so that \,$\liea' = \liem'_{fl} \oplus \liea'_c$\, holds. \,$\liea'_c$\, is flat,
%and because \,$\liea'$\, is a Cartan subalgebra of \,$\liem'$\,, \,$\liea'_c$\, is a Cartan subalgebra of \,$\liem'_c$\,.
%
%Conversely, suppose that \,$\liea'  = \liem'_{fl} \oplus \liea'_c$\, holds with some Cartan subalgebra \,$\liea'_c$\, of \,$\liem'_c$\,.
%Then \,$\liea'$\, is flat because \,$\liem'_{fl}$\, and \,$\liea'_c$\, are flat and \eqref{eq:roots:flat:a'liem'fl} holds, and the maximality of \,$\liea'_c$\,
%implies the maximality of \,$\liea'$\,. 
%
%The statement on \,$\rk(\liem')$\, now follows immediately.
%\beweisende

We now fix a Cartan subalgebra \,$\liea' \subset \liem'$\,; 
%by Proposition~\ref{P:roots:flat} 
because of Equation~\eqref{eq:roots:mm'0} there exists a Cartan subalgebra \,$\liea'_c$\, of \,$\liem'_c$\,
so that 
\begin{equation}
\label{eq:roots:a'}
\liea' = \liem'_{fl} \oplus \liea'_c
\end{equation}
holds.
\,$\liem'_c$\, is of compact type, and therefore we have the usual root space decomposition with respect to \,$\liea'_c$\,: 
We put for any \,$\alpha_c \in (\liea'_c)^*$\,
$$ (\liem'_c)_{\alpha_c} := \Menge{\,X \in \liem'_c\,}{\,\forall Z \in \liea'_c: \ad(Z)^2 X = -\alpha_c(Z)^2X\,} $$
and consider the root system
$$ \Delta(\liem'_c,\liea'_c) := \Menge{\,\alpha_c \in (\liea'_c)^*\setminus \{0\}\,}{\,(\liem'_c)_{\alpha_c} \neq \{0\}\,} \; , $$
then we have
\begin{equation}
\label{eq:roots:mc'decomp}
\liem'_c = \liea'_c \;\oplus\; \bigoplus_{\alpha \in \Delta_+(\liem'_c,\liea'_c)} (\liem'_c)_{\alpha_c} \;,
\end{equation}
where \,$\Delta_+(\liem'_c,\liea'_c)$\, is any system of positive roots in \,$\Delta(\liem'_c,\liea'_c)$\,.

Now we define for any \,$\alpha_c \in (\liea'_c)^*$\, the linear form \,$\alpha \in (\liea')^*$\, by \,$\alpha|\liea'_c = \alpha_c$\, and
\,$\alpha|\liem'_{fl} = 0$\,. Then for any \,$\alpha_c \in (\liea'_c)^* \setminus \{0\}$\,
$$ \liem_\alpha' := \Menge{\,X \in \liem'\,}{\,\forall Z \in \liea' : \ad(Z)^2 X = -\alpha(Z)^2 X \,} \;=\; (\liem'_c)_{\alpha_c} $$
holds, and therefore Equations~\eqref{eq:roots:m'flc}, \eqref{eq:roots:a'} and \eqref{eq:roots:mc'decomp} 
show that we have the root space decomposition
\begin{equation}
\label{eq:roots:m'decomp}
\liem' = \liea' \;\oplus\; \bigoplus_{\alpha \in \Delta_+(\liem',\liea')} \liem'_{\alpha}
\end{equation}
with respect to the root system
$$ \Delta(\liem',\liea') := \Menge{\,\alpha\,}{\,\alpha_c \in \Delta(\liem'_c,\liea'_c)\,} \;; $$
\,$\Delta_+(\liem',\liea')$\, is a system of positive roots in \,$\Delta(\liem',\liea')$\,. 

The following proposition describes relations between the root systems \,$\Delta(\liem,\liea)$\, and \,$\Delta(\liem',\liea')$\,, as well as between
the root spaces \,$\liem_\alpha'$\, and \,$\liem_\lambda$\,. 

\begin{Prop}
\label{P:cla:subroots:subroots-neu}
Let \,$\liea$\, be a Cartan subalgebra of \,$\liem$\, such that \,$\liea' := \liea \cap \liem'$\, is a Cartan subalgebra of \,$\liem'$\,.
\begin{enumerate}
\item
The roots resp.~root spaces of \,$\liem$\, and of \,$\liem'$\, are related by
the following equations:
\begin{gather}
\label{eq:cla:subroots:subroots-neu:liea'}
\liem_0' = \liea' \;, \\
\label{eq:cla:subroots:subroots-neu:toshow-Delta}
\Delta(\liem',\liea') \subset \Menge{\,\lambda|\liea'\,}{\,\lambda \in \Delta(\liem,\liea), \lambda|\liea' \neq 0\,} \;, \\
\label{eq:cla:subroots:subroots-neu:toshow-liemalpha}
% \forall \alpha \in \Delta(\liem',\liea'): \liem_\alpha' \subset \bigoplus_{\substack{\lambda \in \Delta(\liem,\liea) \\ \lambda|\liea' = \alpha}} \liem_\lambda \\
\textstyle
\forall \alpha \in \Delta(\liem',\liea')\;:\; \liem_\alpha' = \left( \bigoplus_{\substack{\lambda \in \Delta(\liem,\liea) \\ \lambda|\liea' = \alpha}} \liem_\lambda \right) \;\cap\; \liem' \; .
\end{gather}
\item
We have \,$\rk(\liem') = \rk(\liem)$\, if and only if \,$\liea' = \liea$\, holds. If this is the case, then we have
\begin{equation}
\label{eq:cla:subroots:subroots-neu:c}
\Delta(\liem',\liea') \subset \Delta(\liem,\liea) \;,\quad \forall \alpha \in \Delta(\liem',\liea') : \liem'_\alpha = \liem_\alpha \cap \liem' 
\end{equation}
\end{enumerate}
\end{Prop}

\beweis
Equation~\eqref{eq:cla:subroots:subroots-neu:liea'} follows by the usual argument: Because \,$\liea'$\, is flat, we have \,$\liea' \subset \liem_0'$\,.
Conversely, let \,$X \in \liem_0'$\, be given. Then we have for every \,$Z \in \liea'$\,: \,$\ad(Z)^2 X = 0$\, and therefore 
\,$0 = \vkap(\ad(Z)^2 X,X) = -\vkap(\ad(Z)X,\ad(Z)X)$\,. 
Because the Killing form \,$\vkap$\, of \,$\lieg$\, is negative definite, it follows that 
\,$\ad(Z)X = 0$\, holds. From this fact and \,$[\liea',\liea']=\{0\}$\,, we see that
\,$[\liea'+\R X,\liea'+\R X] = \{0\}$\, holds, showing that \,$\liea' + \R X$\, is flat. Because of the maximality of \,$\liea'$\,
we conclude \,$X \in \liea'$\,.

Now let \,$\Delta_+ \subset \Delta(\liem,\liea)$\, be a system of positive roots of \,$\liem$\, and put \,$\widetilde{\Delta}_+ := \Delta_+ \cup \{0\}$\,.
Then we have by \eqref{eq:roots:mdecomp} and \eqref{eq:roots:m0a} 
\begin{equation}
\label{eq:cla:subroots:subroots:mdecomp}
\liem = \bigoplus_{\lambda \in \widetilde{\Delta}_+} \liem_\lambda \; . 
\end{equation}

Let \,$\alpha \in \Delta(\liem',\liea')$\, and \,$X \in \liem'_\alpha$\, be given.
We have \,$X \in \liem$\, and therefore Equation~\eqref{eq:cla:subroots:subroots:mdecomp} shows that there exists a decomposition
\begin{equation}
\label{eq:cla:subroots:subroots:wdecomp}
X = \sum_{\lambda \in \widetilde{\Delta}_+} X_\lambda
\end{equation}
with suitable (unique) \,$X_\lambda \in \liem_\lambda$\, for \,$\lambda \in \wt{\Delta}_+$\,.

Because of \,$X \in \liem'_\alpha$\, and Equation~\eqref{eq:cla:subroots:subroots:wdecomp}, we have
\begin{equation}
\label{eq:cla:subroots:subroots:ad1}
\forall Z \in \liea' \; : \;
\ad(Z)^2 X
= -\alpha(Z)^2 X
= -\sum_{\lambda \in \widetilde{\Delta}_+} \alpha(Z)^2 X_\lambda \; .
\end{equation}
On the other hand, we have \,$X_\lambda \in \liem_\lambda$\, for every \,$\lambda \in \widetilde{\Delta}_+$\, and therefore
\begin{equation}
\label{eq:cla:subroots:subroots:ad2}
\forall Z \in \liea' \; : \;
\ad(Z)^2 X = \sum_{\lambda \in \widetilde{\Delta}_+} \ad(Z)^2 X_\lambda
= -\sum_{\lambda \in \widetilde{\Delta}_+} \lambda(Z)^2 X_\lambda \; .
\end{equation}
By comparing Equations~\eqref{eq:cla:subroots:subroots:ad1} and \eqref{eq:cla:subroots:subroots:ad2} we obtain
$$ \forall Z \in \liea' \; : \; \sum_{\lambda \in \widetilde{\Delta}_+} \alpha(Z)^2 X_\lambda = \sum_{\lambda \in \widetilde{\Delta}_+} \lambda(Z)^2 X_\lambda  $$
and therefore, because of the directness of the sum in Equation~\eqref{eq:cla:subroots:subroots:mdecomp},
$$ \forall \lambda \in \widetilde{\Delta}_+,\;Z \in \liea' \; : \; \alpha(Z)^2 \cdot X_\lambda = \lambda(Z)^2 \cdot X_\lambda \;.$$
Thus, we have \,$X_\lambda = 0$\, for every \,$\lambda \in \widetilde{\Delta}_+$\,
with \,$\,\lambda^2|\liea' \neq \alpha^2$\,, and therefore Equation~\eqref{eq:cla:subroots:subroots:wdecomp} shows
$$ X
\in \bigoplus_{\substack{\lambda \in \widetilde{\Delta}_+ \\ \lambda^2|\liea' = \alpha^2}} \liem_\lambda
= \bigoplus_{\substack{\lambda \in \Delta(\liem,\liea) \\ \lambda|\liea' = \alpha}} \liem_\lambda \; ;
$$
for the last equality, one has to note that for any pair \,$(\alpha_1, \alpha_2)$\, of linear forms
on \,$\liea'$\,, \,$\alpha_1^2 = \alpha_2^2$\, already implies \,$\alpha_1 = \pm \alpha_2$\,. This completes the proof of 
the inclusion ``$\subset$'' of Equation \eqref{eq:cla:subroots:subroots-neu:toshow-liemalpha}; its converse inclusion
follows immediately from the definitions of \,$\liem_\lambda$\, and \,$\liem_\alpha'$\,. 

For any given \,$\alpha \in \Delta(\liem',\liea')$\, we have \,$\liem_\alpha' \neq \{0\}$\, and therefore
\eqref{eq:cla:subroots:subroots-neu:toshow-liemalpha} implies the existence of \,$\lambda \in \Delta(\liem,\liea)$\, with \,$\alpha = \lambda|\liea'$\,;
this observation proves \eqref{eq:cla:subroots:subroots-neu:toshow-Delta}.

It remains to verify (b). We suppose \,$\rk(\liem') = \rk(\liem)$\,. Because \,$\liea'$\, and \,$\liea$\, are Cartan subalgebras
of \,$\liem'$\, and \,$\liem$\, respectively, we then have
\,$\dim \liea'  = \rk(\liem') = \rk(\liem) = \dim \liea$\, and therefore \,$\liea' = \liea$\,.
The remaining statements of (b) now follow from (a).
\beweisende

\begin{Def}
\label{D:cla:subroots:Elemcomp}
Let \,$\liea$\, be a Cartan subalgebra of \,$\liem$\,
so that \,$\liea' := \liea \cap \liem'$\, is a Cartan subalgebra of \,$\liem'$\,. Also let \,$\alpha \in \Delta(\liem',\liea')$\, be given.
Recall that by Proposition~\ref{P:cla:subroots:subroots-neu}(a) there exists at least one root \,$\lambda \in \Delta(\liem,\liea)$\, with
\,$\lambda|\liea' = \alpha$\,. We call \,$\alpha$\, 
\begin{enumerate}
\item \emph{elementary}, if there is only one root \,$\lambda \in \Delta(\liem,\liea)$\, with \,$\lambda|\liea' = \alpha$\,;
\item \emph{composite}, if there are at least two different roots \,$\lambda, \mu \in \Delta(\liem,\liea)$\, with \,$\lambda|\liea' = \mu|\liea' = \alpha$\,.
\end{enumerate}
\end{Def}

In the situation described in Definition~\ref{D:cla:subroots:Elemcomp}, elementary roots play a special role:
If \,$\alpha \in \Delta(\liem',\liea')$\, is elementary, then the root space \,$\liem_\alpha'$\, 
is contained in the root space \,$\liem_\lambda$\,, where \,$\lambda \in \Delta(\liem,\liea)$\, is the unique root with \,$\lambda|\liea' = \alpha$\,.
% (see Proposition~\ref{P:cla:subroots:subroots-neu}(a)). 
As we will see in Proposition~\ref{P:cla:subroots:Comp} below, this property
causes restrictions for the possible positions (in relation to \,$\liea'$\,) of \,$\lambda$\,. The exploitation of these restrictions will play
an important role in the classification of the Lie triple systems of the complex quadric in Section~\ref{Se:cla}.

It should be mentioned that in the case \,$\rk(\liem') = \rk(M)$\, we have \,$\liea' = \liea$\,, and therefore in that case
every \,$\alpha \in \Delta(\liem',\liea')$\, is elementary (see Proposition~\ref{P:cla:subroots:subroots-neu}(b)).

For any linear form \,$\lambda \in \liea^*$\, we now denote by \,$\lambda^\sharp$\, the Riesz vector
corresponding to \,$\lambda$\,, i.e.~the vector \,$\lambda^\sharp \in \liea$\, characterized by \,$\g{\,\cdot\,}{\lambda^\sharp} = \lambda$\,. 

\begin{Prop}
\label{P:cla:subroots:Comp}
Let \,$\liea$\, be a Cartan subalgebra of \,$\liem$\,
so that \,$\liea' := \liea \cap \liem'$\, is a Cartan subalgebra of \,$\liem'$\,. Also let \,$\alpha \in \Delta(\liem',\liea')$\, 
be given. 
\begin{enumerate}
\item
If \,$\alpha$\, is elementary 
and \,$\lambda \in \Delta(\liem,\liea)$\, is the unique root with \,$\lambda|\liea' = \alpha$\,, then we have \,$ \lambda^\sharp \in \liea'$\,. 
\item
If \,$\alpha$\, is composite and \,$\lambda, \mu \in \Delta(\liem,\liea)$\, are two different roots with
\,$\lambda|\liea' = \alpha = \mu|\liea'$\,, then \,$\lambda^\sharp - \mu^\sharp$\, is orthogonal to \,$\liea'$\,. 
\end{enumerate}
\end{Prop}

\beweis
\emph{For (a).}
Let \,$\alpha \in \Delta(\liem',\liea')$\, be an elementary root and \,$\lambda \in \Delta(\liem,\liea)$\, be the root
with \,$\lambda|\liea' = \alpha$\,. 
We fix \,$X \in \Sph(\liem_\alpha')$\, arbitrarily.%
\footnote{For any euclidean space \,$V$\,, we denote the unit sphere in \,$V$\, by \,$\Sph(V)$\,.}
By Proposition~\ref{P:cla:subroots:subroots-neu}(a) we have \,$X \in \liem_\lambda$\,.
Thus there is exactly one \,$\widehat{X} \in \liek_\lambda$\, which is related to \,$X \in \liem_\lambda$\, in the sense (see \cite{Loos:1969-2}, Lemma~VI.1.5, p.~62) that 
\begin{equation}
\label{eq:cla:subroots:Comp:what1}
\forall Z \in \liea \; : \; \bigr( \; [Z,X] = \lambda(Z)\cdot \wh{X} \qmq{and} [Z,\wh{X}] = -\lambda(Z)\cdot X \;\bigr) 
\end{equation}
holds, and then we also have
\begin{equation}
\label{eq:cla:subroots:Comp:what2}
[X,\wh{X}] = \lambda^\sharp \; . 
\end{equation}
We now fix \,$Z \in \liea'$\, so that \,$\lambda(Z) = \alpha(Z) = -1$\, holds; then we have
$$ \liem' \;\overset{(*)}{\ni}\; \ad(X)^2 Z  = [X,[X,Z]] = - [X,[Z,X]] \overset{\eqref{eq:cla:subroots:Comp:what1}}{=} [X,\wh{X}]
\overset{\eqref{eq:cla:subroots:Comp:what2}}{=} \lambda^\sharp \; , $$ 
where $(*)$ follows from the fact that \,$\liem'$\, is a Lie triple system. Therefore we have \,$\lambda^\sharp \in \liem' \cap \liea = \liea'$\,. 

(b) is obvious.
\beweisende

\begin{Remark}
Investigating root systems of Lie algebras, \textsc{Eschenburg} used concepts similar to the elementary/composite roots from Definition~\ref{D:cla:subroots:Elemcomp}, 
see \cite{Eschenburg:1984}, Abschnitt 91, p.~131ff.\,.
That situation is different from ours, because in contrary to symmetric spaces, the root spaces of Lie algebras are always \,$1$-dimensional.
\end{Remark}

We also consider the Weyl group of a Lie triple system:

\begin{Def}
\label{D:roots:weyl}
Let \,$\liea'$\, be a Cartan subalgebra of \,$\liem'$\,. 
For \,$\alpha \in \Delta(\liem',\liea')$\, we denote by \,$R_\alpha' : \liea' \to \liea'$\, the orthogonal reflection in the hyperplane
\,$\alpha^{-1}(\{0\})$\,. Then we call the group of orthogonal transformations of \,$\liea'$\, generated by \,$\Menge{R_\alpha'}{\alpha \in \Delta(\liem',\liea')}$\,
the \emph{Weyl group} \,$W(\liem',\liea')$\, of \,$\liem'$\, (with respect to \,$\liea'$\,). \,$W(\liem',\liea')$\, also acts on \,$(\liea')^*$\,
via the action \,$(B,\alpha) \mapsto \alpha \circ B^{-1}$\,. 

For \,$\liem$\, we use the analogous notations, where we omit the \,$'$\, from the symbols.
\end{Def}

The following proposition shows that various well-known facts concerning the Weyl group of a Riemannian symmetric space of compact type
transfer to the present situation for \,$\liem'$\,. It also gives a relation between the Weyl groups \,$W(\liem',\liea')$\, and \,$W(\liem,\liea)$\,.

\newpage

\begin{Prop}
\label{P:roots:weyl}
Let \,$\liea$\, be a Cartan subalgebra of \,$\liem$\, so that \,$\liea' := \liea \cap \liem'$\, is a Cartan subalgebra of \,$\liem'$\,. 
\begin{enumerate}
\item
We have \,$W(\liem',\liea') \subset \Menge{\;(P_{\liea'} \circ B)|\liea'\;}{\;B \in W(\liem,\liea)\;}$\,, where \,$P_{\liea'}: \liea \to \liea$\, denotes
the orthogonal projection onto \,$\liea'$\,. 

In particular, in the case \,$\rk(\liem') = \rk(\liem)$\, we have \,$W(\liem',\liea) \subset W(\liem,\liea)$\, (remember that \,$\liea' = \liea$\, then holds).
\item
Let us denote by \,$K'$\, the connected Lie subgroup of \,$G$\, with Lie algebra \,$\liek' := [\liem',\liem']$\,, then \,$K'$\, is also a subgroup of \,$K$\,. 
Also let \,$B \in W(\liem',\liea')$\, be given. Then there exists \,$g \in K'$\, so that \,$B = \Ad(g)|\liea'$\, holds.
\item
The Weyl group \,$W(\liem',\liea')$\, leaves the root system \,$\Delta(\liem',\liea')$\, invariant.
\end{enumerate}
\end{Prop}

\beweis
\emph{For (a).}
It suffices to show \,$R_\alpha' = (P_{\liea'} \circ R_\lambda)|\liea'$\, for any given \,$\alpha \in \Delta(\liem',\liea')$\,, where
\,$\lambda \in \Delta(\liem,\liea)$\, is chosen so that \,$\lambda|\liea' = \alpha$\, holds; the existence of such a \,$\lambda$\, follows from 
Equation~\eqref{eq:cla:subroots:subroots-neu:toshow-Delta} in Proposition~\ref{P:cla:subroots:subroots-neu}(a). Indeed, in this situation we
have \,$\alpha^\sharp = P_{\liea'}(\lambda^\sharp)$\, and therefore for every \,$Z \in \liea'$\,
$$ R_\alpha'(Z) = Z - 2\g{Z}{\alpha^\sharp}\alpha^\sharp = Z - 2\g{Z}{P_{\liea'}(\lambda^\sharp)}P_{\liea'}(\lambda^\sharp) 
\overset{(*)}{=} P_{\liea'}(\,Z - 2\g{Z}{\lambda^\sharp}\lambda^\sharp\,) = P_{\liea'}(R_\lambda(Z)) \; ; $$
the equals sign marked $(*)$ follows from the fact that \,$Z \in \liea'$\, holds.

\emph{For (b).}
It suffices to show that for any given \,$\alpha \in \Delta(\liem',\liea')$\,, there exists \,$g \in K'$\, so that \,$R_\alpha' = \Ad(g)|\liea'$\, holds.
For this purpose we choose \,$Z_0 \in \liea'$\, so that \,$\alpha(Z_0) = 1$\, holds, moreover 
we fix \,$X \in \Sph(\liem_\alpha')$\, and put
$$ \wh{X} := [Z_0,X] \in \liek' $$
(then \,$\wh{X}$\, is related to \,$X$\, in the sense of \cite{Loos:1969-2}, Lemma~VI.1.5, p.~62) and 
\begin{equation}
\label{eq:roots:weyl:g}
g := \Exp(t_0\,X) \in K' \qmq{with} t_0 := \frac{\pi}{\|\alpha^\sharp\|} \;,
\end{equation}
where \,$\Exp: \lieg \to G$\, denotes the exponential map of the Lie group \,$G$\,. We then have \,$\Ad(g)\alpha^\sharp = -\alpha^\sharp$\,
(compare \cite{Loos:1969-2}, Lemma~VI.1.5(c), p.~62) and by a calculation similar to the one in the proof of \cite{Loos:1969-2}, Lemma~VI.1.5(c) on p.~63
for every \,$Z \in \liea'$\, which is orthogonal to \,$\alpha^\sharp$\,: \,$\Ad(g)Z=Z$\,. Therefrom \,$\Ad(g)|\liea' = R_\alpha'$\, follows.

\emph{For (c).}
It suffices to show that we have \,$\beta \circ B^{-1} \in \Delta(\liem',\liea')$\, for any given \,$\beta \in \Delta(\liem',\liea')$\, and \,$B \in W(\liem',\liea')$\,.
By (b), there exists \,$g \in K'$\, so that \,$B = \Ad(g)|\liea'$\, holds. In particular, we have \,$\Ad(g)\liea' = \liea'$\,, and below we will show 
\begin{equation}
\label{eq:roots:weyl:Adgm'}
\Ad(g)\liem' = \liem' \; . 
\end{equation}
Therefrom we obtain
\begin{align*}
\Ad(g)\liem_{\beta}' & = \Menge{\Ad(g)X}{X \in \liem',\;\forall Z \in \liea': \ad(Z)^2X = -\beta(Z)^2\,X} \\
& = \Menge{X \in \Ad(g)\liem'}{\forall Z \in \liea': \ad(Z)^2(\Ad(g)^{-1}X) = -\beta(Z)^2\,(\Ad(g)^{-1}X)} \\
& = \Menge{X \in \liem'}{\forall Z \in \liea': \Ad(g)^{-1}\bigr(\;\ad(\Ad(g)Z)^2X \;\bigr) = \Ad(g)^{-1}\bigr(\;-\beta(Z)^2 X\;\bigr)} \\
& = \Menge{X \in \liem'}{\forall Z \in \liea': \ad(\Ad(g)Z)^2X = -\beta(Z)^2 X} \\
& = \Menge{X \in \liem'}{\forall Z \in \Ad(g)\liea': \ad(Z)^2X = -(\beta\circ \Ad(g)^{-1})(Z)^2 X} \\
& = \Menge{X \in \liem'}{\forall Z \in \liea': \ad(Z)^2X = -(\beta\circ B^{-1})(Z)^2 X} = \liem_{\beta \circ B^{-1}} \; .
\end{align*}
Because of \,$\beta \in \Delta(\liem',\liea')$\, we have \,$\liem_\beta' \neq \{0\}$\, and therefore by the preceding calculation also \,$\liem_{\beta \circ B^{-1}}
\neq \{0\}$\,, whence \,$\beta \circ B^{-1} \in \Delta(\liem',\liea')$\, follows.

For the proof of Equation~\eqref{eq:roots:weyl:Adgm'}, we may suppose without loss of generality that \,$B = R_\alpha'$\, holds for some
\,$\alpha \in \Delta(\liem',\liea')$\,. We now use the notations from the proof of (b), especially \,$g$\, is now given by \eqref{eq:roots:weyl:g},
let \,$Y \in \liem'$\, be given and consider the function
$$ f: \R \to \liem,\; t \mapsto \Ad(\Exp(t\,\wh{X}))Y = \exp(t\,\ad(\wh{X}))Y \; , $$
where \,$\exp: \End(\lieg) \to \GL(\lieg)$\, is the usual exponential map of endomorphisms. \,$f$\, solves the differential equation
\begin{equation}
\label{eq:cla:subroots:subroots-neu:DGL}
y' = \ad(\wh{X}) y \; .
\end{equation}
Because \,$\liem'$\, is a Lie triple system, and we have \,$\wh{X} \in \liek' = [\liem',\liem']$\,, we see that the endomorphism \,$\ad(\wh{X})$\, leaves
\,$\liem'$\, invariant. Because we also have \,$f(0) = Y \in \liem'$\,, the solution \,$f$\, of the differential
equation~\eqref{eq:cla:subroots:subroots-neu:DGL} runs entirely in \,$\liem'$\,. In particular we have \,$\Ad(g)Y = f(t_0) \in \liem'$\,.
Thus we have shown \,$\Ad(g)\liem' \subset \liem'$\,; because \,$\Ad(g)$\, is a linear isomorphism, we conclude \eqref{eq:roots:weyl:Adgm'}.
\beweisende

\section{The geometry of the complex quadric}
\label{Se:Q}

We now turn our attention specifically to the complex quadric, which we also regard as a complex hypersurface in the complex projective
space \,$\CP^{m+1}$\,:
$$ Q := Q^m := \left\{ \;[z_0,\dotsc,z_{m+1}] \in \CP^{m+1} \;\left| \; \sum_k z_k^2 = 0 \; \right. \right\} \; . $$
We regard \,$\CP^{m+1}$\, as a Hermitian manifold via the Fubini-Study metric \,$\g{\cdot}{\cdot}$\, and the usual complex structure \,$J$\,.
These data are characterized by the fact that the Hopf fibration \,$\pi: \Sph^{2m+3} \to \CP^{m+1}, \; z \mapsto [z] := \C z$\, 
(where we regard the unit sphere \,$\Sph^{2m+3}$\, as a submanifold of \,$\C^{m+2}$\,) becomes a Hermitian submersion, meaning that if we consider for
\,$z \in \Sph(\C^{m+2})$\, the horizontal space \,$\scrH_z := \ker(T_z \pi)^{\perp}$\,, then \,$T_z\pi|\scrH_z: \scrH_z \to T_{\pi(z)}\CP^{m+1}$\,
is a complex linear isometry. It should be noted that via the mentioned structures, we also obtain
complex inner products \,$\gC{\cdot}{\cdot}$\, on the tangent spaces of \,$\CP^{m+1}$\,:
$$ \forall p \in \CP^{m+1},\; v,w \in T_p\CP^{m+1} \; : \; \gC{v}{w} = \g{v}{w} + i\cdot \g{v}{Jw} \; . $$

One objective of the present paper is to classify the totally geodesic submanifolds of \,$Q$\,. To this end, we require some information
about the geometry of \,$Q$\,, which is now given. Most of the results stated in the present section are taken from
the paper \cite{Reckziegel:quadrik-1995} by \textsc{H.~Reckziegel}; in particular the concept of a \CQ-structure and its application
to the study of the complex quadric were introduced there.

The following proposition describes the fundamental data concerning the intrinsic and extrinsic geometry of \,$Q \subset \CP^{m+1}$\,. 
For \,$p \in Q$\, we denote by \,$\nsp{Q}{p} := (T_pQ)^{\perp,T_p\CP^{m+1}}$\, the normal space of \,$Q \subset \CP^{m+1}$\, at \,$p$\,
and also consider the set of unit normal vectors \,$\unsp{Q}{p} := \Sph(\nsp{Q}{p})$\,. \,$\unsp{Q}{p}$\, is a ``circle'' in the sense
that for \,$\eta_0 \in \unsp{Q}{p}$\,, \,$\unsp{Q}{p} = \Menge{\lambda\eta_0}{\lambda \in \Sph^1}$\, holds.%
\footnote{Here and in the sequel, we regard the unit circle \,$\Sph^1$\, also as a subset of \,$\C$\,.}

\begin{Prop}
\label{P:Q:AR}
Let \,$p \in Q$\, and \,$\eta \in \unsp{Q}{p}$\, be given.
\begin{enumerate}
\item
The shape operator \,$A_\eta: T_pQ \to T_pQ$\, of \,$Q \subset \CP^{m+1}$\, with respect to \,$\eta$\, is an orthogonal, anti-linear%
\footnote{We call an \,$\R$-linear map \,$A: \bbV \to \bbV$\, on the \,$\C$-linear space \,$\bbV$\, \emph{anti-linear}, if \,$A \circ J = -J \circ A$\, holds.}
involution; moreover \,$A_{\lambda\eta} = \lambda\cdot A_\eta$\, holds for \,$\lambda \in \Sph^1$\,. 

An explicit description of \,$A_\eta$\, is given in the following way: For \,$z \in \pi^{-1}\{p\}$\,, we have
\,$\scrH_z = (\C z)^\perp$\,, and the horizontal lift of \,$T_pQ$\, at \,$z$\, is given by
\,$\scrH_zQ := \scrH_z \cap (T_z\pi)^{-1}(T_pQ) = (\C z \operp \C\overline{z})^\perp$\, (where \,$\overline{z}$\, denotes the usual conjugation
of \,$z \in \C^{m+2}$\,). There exists \,$z \in \pi^{-1}(\{p\})$\, (depending on \,$\eta$\,), so that the following diagram commutes:
\begin{equation*}
\begin{minipage}{5cm}
\begin{xy}
\xymatrix{
\scrH_z Q \ar[r]^{v \mapsto \overline{v}} \ar[d]_{\pi_*|\scrH_z Q} & \scrH_z Q \ar[d]^{\pi_*|\scrH_z Q} \\
T_pQ \ar[r]_{A_{\eta}} & T_pQ \; . 
}
\end{xy}
\end{minipage}
\end{equation*}

\item
The curvature tensor \,$R$\, of \,$Q$\, at \,$p$\, is described via \,$A_\eta$\, by the equation
\begin{equation*}
% \forall u,v,w \in T_pQ \; : \; 
R(u,v)w = \gC{w}{v}\,u - \gC{w}{u}\,v - 2\,\g{Ju}{v}\,Jw + \gC{v}{A_\eta w}\,A_\eta u - \gC{u}{A_\eta w}\,A_\eta v \; . 
\end{equation*}
\end{enumerate}
\end{Prop}

\beweis
See \cite{Reckziegel:quadrik-1995}, Section~3.
\beweisende

As Proposition~\ref{P:Q:AR}(b) shows, the curvature tensor of the complex quadric \,$Q$\,, and therefore the local geometric structure of this Riemannian manifold,
is described entirely by the Riemannian metric of \,$Q$\,, its complex structure \,$J$\, and the ``circle of conjugations''
$$ \A(Q,p) := \Menge{A_\eta}{\eta \in \unsp{Q}{p}} = \Menge{\lambda\cdot A_{\eta_0}}{\lambda \in \Sph^1} \;, $$
where \,$\eta_0 \in \unsp{Q}{p}$\, is fixed arbitrarily. In order to perform the investigation of this situation in a clearer manner, we
formulate some relevant concepts and results in an abstract setting.

\begin{Def}
\label{D:Q:CQ}
Let \,$\bbV$\, be a unitary space with complex inner product \,$\gC{\,\cdot\,}{\,\cdot\,}$\, and complex structure \,$J: \bbV \to \bbV, \; v \mapsto i\cdot v$\,.
\begin{enumerate}
\item
A \emph{conjugation} on \,$\bbV$\, is an anti-linear involution \,$A: \bbV \to \bbV$\, which is 
orthogonal with respect to the induced real inner product \,$\g{\,\cdot\,}{\,\cdot\,} := \RE(\gC{\,\cdot\,}{\,\cdot\,})$\,.
\item
Let \,$A$\, be a conjugation on \,$\bbV$\,. Then we call the circle of conjugations \,$\A := \Menge{\lambda A}{\lambda \in \Sph^1}$\,
a \emph{\CQ-structure} on \,$\bbV$\,. We call the pair \,$(\bbV,\A)$\, (or simply \,$\bbV$\,, if the implied \CQ-structure is obvious)
a \emph{\CQ-space}.
\end{enumerate}
We now suppose that \,$(\bbV,\A)$\, and \,$(\bbV',\A')$\, are \CQ-spaces.
\begin{enumerate}
\addtocounter{enumi}{2}
\item
We call a unitary map \,$B: \bbV \to \bbV'$\, a \emph{\CQ-isomorphism}, if \,$B \circ A \circ B^{-1} \in \A'$\, holds for every \,$A \in \A$\,.
In the case \,$(\bbV',\A') = (\bbV,\A)$\, we speak of a \emph{\CQ-automorphism}. We denote the group of \CQ-automorphisms of \,$(\bbV,\A)$\,
by \,$\Aut(\A)$\,. 
\item
We call a complex linear subspace \,$U \subset \bbV$\, a \emph{\CQ-subspace}, if it is invariant under some (and then under every) \,$A \in \A$\,.
\item
We call \,$v \in \bbV$\, \emph{isotropic}, if \,$\gC{v}{Av}=0$\, holds for some (and then for every) \,$A \in \A$\,. We call a (real or complex)
linear subspace \,$U \subset \bbV$\, isotropic, if every \,$v \in U$\, is isotropic.
\end{enumerate}
\end{Def}

\begin{Examples}
\label{E:Q:CQ}
\begin{enumerate}
\item
For \,$p \in Q$\,, \,$(T_pQ, \A(Q,p))$\, is a \CQ-space.
\item
The usual conjugation \,$A_0: \C^n \to \C^n,\;z \mapsto \overline{z}$\, on \,$\C^n$\, is a conjugation in the sense of Definition~\ref{D:Q:CQ}(a) on this
unitary space. Therefore \,$\A_0 := \Menge{\lambda\,A_0}{\lambda \in \Sph^1}$\, is a \CQ-structure on \,$\C^n$\,. We call \,$A_0$\, resp.~\,$\A_0$\,
the \emph{standard conjugation} resp.~the \emph{standard \CQ-structure} on \,$\C^n$\,. 
\end{enumerate}
\end{Examples}

\begin{Prop}
\label{P:Q:CQ}
\begin{enumerate}
\item
Let \,$\bbV$\, be a unitary space and \,$A$\, be a conjugation on \,$\bbV$\,. Then \,$A$\, is 
self-adjoint with respect to \,$\g{\cdot}{\cdot}$\, and therefore real orthogonally diagonalizable. 
If we denote the eigenspace \,$\Eig(A,1)$\, of \,$A$\, by \,$V(A)$\,, we have \,$\Eig(A,-1) = JV(A)$\,
and the orthogonal decomposition of \,$\bbV$\, into totally real subspaces
\begin{equation}
\label{eq:Q:CQ:VAdecomp}
\bbV = V(A) \operp JV(A) \; . 
\end{equation}
For any \,$v \in \bbV$\, we have the decomposition
\begin{equation}
\label{eq:Q:CQ:vdecomp}
v = \RE_A(v) + J\,\IM_A(v) \qmq{with} \RE_A(v) := \tfrac{1}{2}(Av+v)\in V(A) \qmq{and} \IM_A(v) := \tfrac{1}{2}J(Av-v) \in V(A) \; .
\end{equation}
\item
A linear subspace \,$U$\, of a \CQ-space \,$(\bbV,\A)$\, is a complex-$k$-dimensional \CQ-subspace if and only if there exists \,$A \in \A$\, and a
real-$k$-dimensional subspace \,$W \subset V(A)$\, so that \,$U = W \operp JW$\, holds. If \,$U$\, is a \CQ-subspace, then this representation
can be achieved for every \,$A \in \A$\,.
\item
Let \,$(\bbV,\A)$\, be a \CQ-space.
\begin{enumerate}
\item 
A complex linear subspace \,$U \subset \bbV$\, of complex dimension \,$k$\, is isotropic if and only if for some (and then for every) \,$A \in \A$\,,
there exists a real-$2k$-dimensional linear subspace \,$W \subset V(A)$\, and an orthogonal complex structure \,$\tau: W \to W$\, so that
\,$U = \Menge{x + J(\tau x)}{x \in W}$\, holds. 
\item
A totally real linear subspace \,$U \subset \bbV$\, of real dimension \,$k$\, is isotropic if and only if for some (and then for every) \,$A \in \A$\,,
there exist real-$k$-dimensional linear subspaces \,$W_1,W_2 \subset V(A)$\, with \,$W_1 \perp W_2$\, and an \,$\R$-linear isometry
\,$\tau: W_1 \to W_2$\, so that \,$U = \Menge{x + J(\tau x)}{x \in W_1}$\, holds. 
\item
If \,$U \subset \bbV$\, is an isotropic subspace, then its ``complex closure'' \,$U + JU$\, is also isotropic.
\end{enumerate}
\end{enumerate}
\end{Prop}

\beweis
\emph{For (a).} For any \,$v,w \in \bbV$\, we have \,$\g{Av}{w} = \g{Av}{A^2 w} = \g{v}{Aw}$\,, thus \,$A$\, is self-adjoint with respect to
\,$\g{\cdot}{\cdot}$\,, and therefore real orthogonally diagonalizable. Because \,$A$\, is involutive, \,$1$\, and \,$-1$\, are its only
possible eigenvalues, and because \,$A$\, is anti-linear, we have \,$\Eig(A,-1) = J(\Eig(A,1))$\,. This shows that \,$V(A) := \Eig(A,1)$\, is
a totally real subspace of \,$\bbV$\, and that the decomposition~\eqref{eq:Q:CQ:VAdecomp} holds. For \,$v \in \bbV$\,, it is obvious that
\,$\RE_A(v), \IM_A(v) \in V(A)$\, holds (with \,$\RE_A(v),\IM_A(v)$\, defined as in \eqref{eq:Q:CQ:vdecomp}), and that we have \,$v = \RE_A(v) + J\,\IM_A(v)$\,.

\emph{For (b).} 
Suppose that \,$U$\, is a \,$k$-dimensional \CQ-subspace of \,$\bbV$\,, let \,$A \in \A$\, be given and put \,$W := U \cap V(A)$\,. We will show
that \,$U = W \operp JW$\, holds with this choice of \,$W$\,; from this equation it also follows that \,$W$\, is of real dimension \,$k$\,.
\,$W \operp JW \subset U$\, holds simply because of \,$W \subset U$\, and \,$U$\, is a complex linear subspace.
For the converse inclusion, let \,$v \in U$\, be given. Because \,$U$\, is a \CQ-subspace
of \,$\bbV$\,, we have \,$x := \RE_A v \in U \cap V(A) = W$\, and similarly \,$y := \IM_A v \in W$\,. This shows that
\,$v = x+Jy \in W \operp JW$\, holds.

Conversely, if \,$U$\, is a linear subspace of \,$\bbV$\, so that \,$U = W \operp JW$\, holds with some linear subspace \,$W \subset V(A)$\,
and \,$A \in \A$\,, then \,$U$\, is clearly complex and \,$A$-invariant, and therefore also invariant under every other \,$A' \in \A$\,.
Hence, \,$U$\, is a \CQ-subspace of \,$\bbV$\,.

\emph{For (c).}
At first, we let \,$U \subset \bbV$\, be any (complex or real linear) isotropic subspace. Then for \,$A \in \A$\,, \,$\beta: U \times U \to \C, \;
(v,w) \mapsto \gC{v}{Aw}$\, is a symmetric, \,$\R$-bilinear map, whose corresponding quadratic form vanishes because \,$U$\, is isotropic.
Therefore we have for any \,$v,w \in U$\,: \,$\beta(v,w) = \tfrac{1}{2}(\beta(v+w,v+w) - \beta(v,v) - \beta(w,w)) = 0$\, and hence
\begin{equation}
\label{eq:Q:CQ:iso-vAw}
\forall A\in \A,\; v,w \in U \; : \; \gC{v}{Aw} = 0 \; . 
\end{equation}
Let us now consider \,$\wh{U} := U+JU$\, and 
let \,$\widehat{v} \in \widehat{U}$\, be given, say \,$\widehat{v} = v_1 + Jv_2$\, with \,$v_1,v_2 \in U$\,. Then we have
\begin{align}
\gC{\widehat{v}}{A\widehat{v}}
& = \gC{v_1 + Jv_2}{Av_1 - JAv_2}
= \gC{v_1}{Av_1} - \gC{v_1}{JAv_2} + \gC{Jv_2}{Av_1} - \gC{Jv_2}{JAv_2} \notag \\
\label{eq:Q:CQ:iso-Udach}
& = \gC{v_1}{Av_1} - \gC{v_1}{JAv_2} - \gC{v_2}{JAv_1} - \gC{v_2}{Av_2} \; .
\end{align}
The first and the fourth summand in \eqref{eq:Q:CQ:iso-Udach} vanish because \,$v_1$\, and \,$v_2$\, are isotropic; the second and the third
summand vanish by \eqref{eq:Q:CQ:iso-vAw}, note that \,$J\circ A = i \, A \in \A$\, holds. Thus we have shown \,$\gC{\widehat{v}}{A\widehat{v}} = 0$\,, and hence
\,$\widehat{U}$\, is isotropic. This proves (c)(iii).

From \eqref{eq:Q:CQ:iso-vAw} one can derive the following equations for every \,$v,w \in U$\, and \,$A \in \A$\, by considering the decomposition
\,$v = \RE_A v + J\IM_A v$\, and the analogous decomposition for \,$w$\,:
\begin{align}
\label{eq:Q:CQ:iso-ReRe}
\g{\RE_A v}{\RE_A w} & = \g{\IM_A v}{\IM_A w} = \tfrac{1}{2}\g{v}{w} \; ,  \\
\label{eq:Q:CQ:iso-ReIm}
\g{\RE_A v}{\IM_A w} & = -\g{\IM_A v}{\RE_A w} \; . 
\end{align}
%Indeed, we have
%\begin{equation}
%\label{eq:CQ:iso:iso:vw}
%\gR{v}{w} = \gR{\RE_A v + J\IM_A v}{\RE_A w + J\IM_Aw} = \gR{\RE_A v}{\RE_A w} + \gR{\IM_A v}{\IM_A w} \; ,
%\end{equation}
%and because of \eqref{eq:Q:CQ:iso-vAw}
%\begin{align*}
%0
%& = \gC{v}{Aw} = \gC{\RE_A v + J\IM_A v}{\RE_A w - J\IM_A w} \\
%& = \gR{\RE_A v}{\RE_A w} - \gR{\IM_A v}{\IM_A w} + i\cdot (\gR{\RE_A v}{\IM_A w} + \gR{\IM_A v}{\RE_A w})
%\end{align*}
%and consequently
%\begin{align}
%\label{eq:CQ:iso:iso:vAw1}
%\gR{\RE_A v}{\RE_A w} & = \gR{\IM_A v}{\IM_A w} \; , \\
%\label{eq:CQ:iso:iso:vAw2}
%\gR{\RE_A v}{\IM_A w} & = - \gR{\IM_A v}{\RE_A w} \; .
%\end{align}
%By combining Equations~\eqref{eq:CQ:iso:iso:vw} and \eqref{eq:CQ:iso:iso:vAw1} we obtain \eqref{eq:Q:CQ:iso-ReRe}, whereas
%Equation~\eqref{eq:CQ:iso:iso:vAw2} proves \eqref{eq:Q:CQ:iso-ReIm}.

\eqref{eq:Q:CQ:iso-ReRe} shows that for \,$v \in U$\, either of the conditions \,$\RE_A v = 0$\, and \,$\IM_A v = 0$\, implies
\,$v = 0$\,. Therefore the surjective \,$\R$-linear maps
$$ \scrR := (\RE_A|U): U \to \RE_A(U) =: W_1 \qmq{and} \scrI := (\IM_A|U): U \to \IM_A(U) =: W_2 $$
are linear isomorphisms, and consequently the linear map \,$\tau := \scrI \circ \scrR^{-1}: W_1 \to W_2$\, also is a linear isomorphism and
$$ U = \Menge{x + J(\tau x)}{x \in W_1} $$
holds; moreover \eqref{eq:Q:CQ:iso-ReRe} shows that \,$\tau$\, is a linear isometry.

If \,$U$\, is a complex isotropic subspace, we let \,$x \in W_1$\, be given. 
Then we have \,$v := x + J(\tau x) \in U$\, and thus also \,$Jv \in U$\,. 
 \,$Jv$\, can therefore be calculated in two different ways:
\begin{align*}
Jv & = J(x + J(\tau x)) = -\tau x + Jx \\
& = \RE_A(Jv) + J\tau(\RE_A(Jv)) \; ;
\end{align*}
and therefore we obtain
\begin{equation}
\label{eq:CQ:iso:iso:cpl}
\RE_A(Jv) = -\tau x \qmq{and} \tau(\RE_A(Jv)) = x \; ,
\end{equation}
whence \,$x = \tau(\RE_A(Jv)) \in W_2$\, follows.
By varying \,$x$\,, we obtain \,$W_1 \subset W_2$\,; because \,$W_1$\, and \,$W_2$\, have the same
dimension (\,$\tau$\, is an isomorphism between them), it follows that we have \,$W_1 = W_2 =: W$\,. Equations~\eqref{eq:CQ:iso:iso:cpl} also shows
that \,$ \tau(\tau x) = -\tau(\RE_A(Jv)) = -x $\, holds for \,$x \in W$\, and therefore \,$\tau$\, is a complex structure on \,$W$\,.
This proves the representation of \,$U$\, given in (c)(i).

On the other hand, if \,$U$\, is a totally real isotropic subspace, we have for any \,$v,w \in U$\,
\begin{align*}
0 & = \g{v}{Jw} = \g{\RE_A v + J\IM_A v}{-\IM_A w + J\RE_A w} \\
& = -\g{\RE_A v}{\IM_A w} + \g{\IM_A v}{\RE_A w} \overset{\eqref{eq:Q:CQ:iso-ReIm}}{=} 2\g{\RE_A w}{\IM_A v} 
\end{align*}
and therefore \,$W_1 \perp W_2$\,. This proves the representation of \,$U$\, given in (c)(ii).

It remains to show the converse implications in the statements of (c)(i) and (c)(ii). For this let \,$A \in \A$\,, linear subspaces \,$W_1,W_2 \subset V(A)$\,
and an \,$\R$-linear isometry \,$\tau: W_1 \to W_2$\, be given, and consider the space
$$ U := \Menge{x + J(\tau x)}{x \in W_1} \; . $$
Then we have for any \,$A \in \A$\, and any \,$v \in U$\,, say \,$v = x + J(\tau x)$\, with \,$x \in W_1$\,
$$ \gC{v}{Av} = \gC{x + J(\tau x)}{x - J(\tau x)} = \g{x}{x} - \g{\tau x}{\tau x} + i\cdot( \g{\tau x}{x} + \g{x}{\tau x} ) = 2i \cdot \g{x}{\tau x} \; . $$
In either of the two cases (i) \,$W_1 = W_2 =: W$\, and \,$\tau: W \to W$\, is a complex structure; (ii) \,$W_1 \perp W_2$\,, we have \,$\g{x}{\tau x} = 0$\,
for all \,$x \in W_1$\,, and therefore \,$U$\, is then isotropic. Moreover, it is easily seen that \,$U$\, is complex in case (i), totally real in case (ii).
\beweisende

We now return to the specific situation of the complex quadric. 
Any \,$B \in \SO(m+2)$\, acts (via complexification) isometrically on \,$\Sph^{2m+3} \subset \C^{m+2}$\,, and thereby 
gives rise to a holomorphic isometry \,$\underline{B}: \CP^{m+1} \to \CP^{m+1}$\, characterized by the equation
\,$\pi \circ (B|\Sph^{2m+3}) = \underline{B} \circ \pi$\, (where \,$\pi: \Sph^{2m+3} \to \CP^{m+1}$\, again is the Hopf fibration). 
\,$B$\, preserves the quadratic form \,$\beta: \C^{m+2}\times \C^{m+2} \to \C,\;(\,(z_k)\,,\,(w_k)\,) \mapsto \sum z_k w_k$\,,
and therefore \,$\underline{B}$\, leaves \,$Q$\, invariant; moreover \,$\underline{B}|Q: Q \to Q$\, is a holomorphic isometry.
In the sequel we regard \,$G := \SO(m+2)$\, as acting on \,$Q$\, via \,$\Psi: G \times Q \to Q,\; (B,p) \mapsto \underline{B}(p)$\,. 
Concerning this action, the following facts hold:

\begin{Prop}
\label{P:Q:Qaction}
\begin{enumerate}
\item \,$G$\, acts transitively on \,$Q$\,. For given \,$p \in Q$\,, say \,$p = \pi(z)$\, with \,$z \in \Sph^{2m+3}$\,,
the isotropy group of this action at \,$p$\, is
$$ K := \Menge{B \in G}{B(W)=W} \;, $$
where \,$W \subset \C^{m+2}$\, is the complex-2-dimensional space spanned by \,$z$\, and \,$\overline{z}$\,.
\item The image of the isotropy representation \,$K \to \Ug(T_pQ),\; B \mapsto T_p\underline{B}$\, is equal to \,$\Aut(\A(Q,p))_0$\,.
\item Let \,$S : \C^{m+2} \to \C^{m+2}$\, the \,$\C$-linear map characterized by \,$S|W = \id_W$\, and \,$S|(W^\perp) = -\id_{W^\perp}$\,. Then we have
\,$-S \in G$\, and the involutive Lie group automorphism
$$ \sigma: G \to G,\; B \mapsto S \circ B \circ S^{-1} $$
satisfies \,$\Fix(\sigma)_0 = K$\,. Consequently \,$\sigma$\, defines the structure of a Hermitian symmetric space on \,$Q$\,;
its canonical covariant derivative is identical to the Levi-Civita covariant derivative of \,$Q$\,.
\item
The Hermitian symmetric space \,$Q$\, is of compact type; it is irreducible for \,$m \neq 2$\,.
\item
The canonical decomposition \,$\lieg = \liek \oplus \liem$\, of the Lie algebra \,$\lieg = \lieo(m+2) \cong \End_-(\C^{m+2})$\, 
of \,$G$\, with respect to \,$\sigma$\, is given by
\begin{align*}
\liek & = \Menge{X \in \lieg}{X(W) \subset W,\;X(W^\perp) \subset W^\perp} \\
\liem & = \Menge{X \in \lieg}{X(W) \subset W^\perp,\;X(W^\perp) \subset W} \; . 
\end{align*}
\end{enumerate}
\end{Prop}

\beweis
For (b), see \cite{Reckziegel:quadrik-1995}, Theorem~2. The remaining facts are well-known.
\beweisende

\begin{Remark}
The isotropy group \,$K$\, mentioned in Proposition~\ref{P:Q:Qaction}(a) is isomorphic to \,$\SO(2) \times \SO(m)$\,. Thus we obtain the
conventional quotient space representation \,$\SO(m+2)/(\SO(2) \times \SO(m))$\, of the complex quadric.
\end{Remark}

Let us now fix \,$p \in Q$\, and consider the corresponding decomposition \,$\lieg = \liek \oplus \liem$\, as in Proposition~\ref{P:Q:Qaction}(d).
As is well-known, the space \,$\liem$\, is linked to the tangent space \,$T_pQ$\, by the linear isomorphism
$$ \tau: \liem \to T_pQ,\; X \mapsto (\Psi(\,\cdot\,,p))_* X_e $$
% $$ \tau: \liem \to T_{p_0}Q,\; X \mapsto \left.\frac{\mathrm{d}\ }{\mathrm{d}t}\right|_{t=0} {\underline{(\Exp(tX))}(p_0)} $$
%is an isomorphism of linear spaces, where \,$\Exp: \lieg \to G$\, is the exponential map of \,$G$\,.
%Remember that for every \,$B \in K$\, (where \,$K$\, is the isotropy group of the action \,$\Psi$\, at \,$p_0$\,),
and we have
\begin{gather}
\label{eq:iso:sym:Ad-tau}
\forall B \in G \; : \; \tau \circ (\Ad(B)|\liem) = T_p \underline{B} \circ \tau \;, \\
\label{eq:iso:sym:R-lie}
\forall X,Y,Z \in \liem \; : \; R(\tau(X),\tau(Y))\tau(Z) = -\tau([[X,Y],Z]) \;, 
\end{gather}
where in the second equation, \,$R$\, again denotes the curvature tensor of \,$Q$\,. 

In the sequel we identify \,$\liem$\, with \,$T_pQ$\, via \,$\tau$\,. In this way, \,$\liem$\, becomes a \CQ-space via the complex structure \,$J$\,,
the complex inner product \,$\gC{\cdot}{\cdot}$\, and the \CQ-structure \,$\A := \A(Q,p)$\,, and we have \,$R(X,Y)Z = -[[X,Y],Z]$\, for any
\,$X,Y,Z \in \liem = T_pQ$\,. Note especially that for a linear subspace \,$\liem' \subset \liem$\,, the notions of \,$\liem'$\, being a Lie triple system
and of \,$\liem'$\, being a curvature-invariant subspace (meaning that \,$R(X,Y)Z \in \liem'$\, holds for every \,$X,Y,Z \in \liem'$\,) therefore coincide.

Via an explicit calculation, one can show that the real inner product \,$\g{\cdot}{\cdot} = \RE(\gC{\cdot}{\cdot})$\, thereby induced
on \,$\liem$\, satisfies
\begin{equation}
\label{eq:Q:skp-killing}
\forall X,Y \in \liem \; : \; \g{X}{Y} = -\tfrac{1}{4m}\cdot \vkap(X,Y) \;,
\end{equation}
where \,$\vkap: \lieg \times \lieg \to \R,\; (X,Y) \mapsto \tr(\ad(X) \circ \ad(Y))$\, is the Killing form of \,$\lieg$\,. Therefore the Riemannian
metric on \,$Q$\, is in accordance with the situation considered in Section~\ref{Se:roots}.

We now suppose \,$m \geq 2$\,.

\begin{Prop}
\label{P:Q:roots}
\begin{enumerate}
\item
The Hermitian symmetric space \,$Q$\, is of rank \,$2$\,, and the Cartan subalgebras \,$\liea$\, of \,$\liem$\,
are exactly the spaces
$$ \liea = \R X \oplus \R JY $$
with \,$A \in \A$\,, \,$X,Y \in \Sph(V(A))$\,, \,$\g{X}{Y} = 0$\,.

\item
Let \,$\liea = \R X \oplus \R J Y$\, be a Cartan subalgebra of \,$\liem$\, as in (a).
Then the following table gives besides \,$\lambda_0 := 0 \in \liea^*$\, a system of positive roots \,$\lambda_k$\, of \,$\liem$\, with respect to \,$\liea$\, (via their Riesz
vectors \,$\lambda_k^\sharp$\,), together with the corresponding root spaces \,$\liem_{\lambda_k}$\, and their multiplicities \,$n_{\lambda_k}$\,:
\begin{center}
\begin{tabular}{|c|c|c|c|}
\hline
\,$k$\, & $\lambda_k^\sharp \in \liea$ & $\liem_{\lambda_k}$ & $n_{\lambda_k}$ \\
\hline
\hline
$0$ & $0$ & $\R X \operp \R JY$ & $2$ \\
$1$ & $\sqrt{2} \cdot JY$ & $J((\R X \operp \R Y)^\perp)$ & $m-2$ \\
$2$ & $\sqrt{2} \cdot X$ & $(\R X \operp \R Y)^\perp$ & $m-2$ \\
$3$ & $\sqrt{2} \cdot (X-JY)$ & $\R(JX+Y)$ & $1$ \\
$4$ & $\sqrt{2} \cdot (X+JY)$ & $\R(JX-Y)$ & $1$ \\
\hline
\end{tabular}
\end{center}
Here \,${}^\perp$\, denotes the ortho-complement in \,$V(A)$\,.
In the case \,$m = 2$\, the roots \,$\lambda_1$\, and \,$\lambda_2$\, do not exist: their multiplicity is zero.
\end{enumerate}
\end{Prop}

\beweis
See \cite{Reckziegel:quadrik-1995}, Sections~5 and 6.
\beweisende

Thus we see that the Riemannian symmetric space \,$Q$\, has (for \,$m \geq 3$\,) the following root diagram:
% (where \,$\lambda_k^\sharp \in \liea$\, denotes the root vector corresponding to the root \,$\lambda_k \in \liea^*$\,):
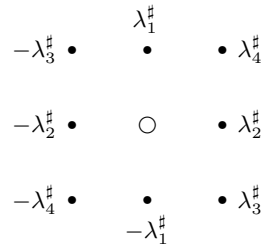
\begin{center}
\strut \\[.3cm]
\setlength{\unitlength}{1cm}
\begin{picture}(2,2)
\put(1,1){\circle{0.2}}
\put(2,1){\circle*{0.1}}       % lambda_2
\put(2,0){\circle*{0.1}}        % lambda_3
\put(1,0){\circle*{0.1}}        % -lambda_1
\put(0,0){\circle*{0.1}}        % -lambda_4
\put(0,1){\circle*{0.1}}        % -lambda_2
\put(0,2){\circle*{0.1}}        % -lambda_3
\put(1,2){\circle*{0.1}}        % lambda_1
\put(2,2){\circle*{0.1}}        % lambda_4
\put(2.2,-0.1){{\small $\lambda_3^\sharp$}}
\put(2.2,0.9){{\small $\lambda_2^\sharp$}}
\put(2.2,1.9){{\small $\lambda_4^\sharp$}}
\put(0.8,2.3){{\small $\lambda_1^\sharp$}}
\put(0.7,-0.5){{\small $-\lambda_1^\sharp$}}
\put(-0.8,-0.1){{\small $-\lambda_4^\sharp$}}
\put(-0.8,0.9){{\small $-\lambda_2^\sharp$}}
\put(-0.8,1.9){{\small $-\lambda_3^\sharp$}}
\end{picture}
\strut \\[1cm]
\end{center}

\section{The classification of Lie triple systems in the complex quadric}
\label{Se:cla}

We continue to consider the specific situation in the complex quadric \,$Q:=Q^m$\, with \,$m \geq 2$\, described in the preceding section. 
Remember that \,$G := \SO(m+2)$\, acts
transitively on \,$Q$\,; we consider the decomposition \,$\lieg = \liek \oplus \liem$\, of the Lie algebra of \,$G$\, induced by the
symmetric structure of \,$Q$\, as was described in Proposition~\ref{P:Q:Qaction} and identify \,$\liem$\, with \,$T_pQ$\, as before.
In this way we again regard \,$\liem$\, as a \CQ-space via the \CQ-structure \,$\A := \A(Q,p)$\,. 

In the present section, we will prove the following theorem:

\begin{Theorem}
\label{T:cla:cla}
A real linear subspace \,$\{0\} \neq \liem' \subsetneq \liem$\, is a Lie triple system if and only if it is of one of the types described in the following list:
%here we have \,$t \in [0,\tfrac{\pi}{4}]$\, and \,$k,k_1,k_2,k'$\, are natural numbers satisfying \,$2 \leq k \leq m$\,, \,$k_1 + k_2 \leq m$\, and
%\,$k' \leq [\tfrac{m}{2}]$\,:
\begin{itemize}
\item[$(\Geo)$]
\,$\liem = \R v$\, holds for some \,$v \in \Sph(\liem)$\,. %  with \,$\vi(v) = t$\,; here we have \,$t \in [0,\tfrac{\pi}{4}]$\,.
\item[$(\Geins,k)$]
\,$\liem'$\, is a \,$k$-dimensional \CQ-subspace of \,$\liem$\, (see Definition~\ref{D:Q:CQ}(d) and Proposition~\ref{P:Q:CQ}(b));
here we have \,$2 \leq k \leq m-1$\,.
\item[$(\Gzwei,k_1,k_2)$]
There exist \,$A \in \A$\, and linear subspaces \,$W_1, W_2 \subset V(A)$\, of real dimension \,$k_1$\, resp.~\,$k_2$\,
so that \,$W_1 \perp W_2$\, and \,$\liem' = W_1 \operp JW_2$\, holds; here we have \,$k_1,k_2 \geq 1$\, and \,$k_1+k_2 \leq m$\,.
\item[$(\Gdrei)$]
There exists \,$A \in \A$\, and an orthonormal system \,$(x,y)$\, in \,$V(A)$\, so that \,$\liem' = \C(x-Jy) \operp \R(x+Jy)$\, holds. %; this type exists only for \,$m \geq 2$\,.
\item[$(\Peins,k)$]
There exists \,$A \in \A$\, so that \,$\liem'$\, is a \,$k$-dimensional \,$\R$-linear subspace of \,$V(A)$\,; here we have \,$1 \leq k \leq m$\,.
\item[$(\Pzwei)$]
There exists \,$A \in \A$\, and \,$x \in \Sph(V(A))$\, so that \,$\liem' = \C x$\, holds.
\item[$(\Atyp)$]
There exists \,$A \in \A$\, and an orthonormal system \,$(x,y,z)$\, in \,$V(A)$\, so that
$$ \liem' = \R(2x+Jy) \operp \R(y+Jx+\sqrt{3}Jz) $$
holds; this type exists only for \,$m \geq 3$\,.
\item[$(\Ieins,k)$]
%There exists \,$A \in \A$\,, a \,$(2k')$-dimensional subspace \,$W \subset V(A)$\, and an orthogonal complex structure \,$j: W \to W$\, so that
%\,$ U = \Menge{x+J(jx)}{x \in W}$\,.
\,$\liem'$\, is a complex \,$k$-dimensional isotropic subspace of \,$\liem$\, (see Definition~\ref{D:Q:CQ}(e) and Proposition~\ref{P:Q:CQ}(c)(i));
here we have \,$1 \leq k \leq \tfrac{m}{2}$\,.
\item[$(\Izwei,k)$]
\,$\liem'$\, is a totally real, real-$k$-dimensional isotropic subspace of \,$\liem$\, (see Definition~\ref{D:Q:CQ}(e) and Proposition~\ref{P:Q:CQ}(c)(ii)); 
here we have \,$1 \leq k \leq \tfrac{m}{2}$\,.
%\,$U$\, is a \,$k'$-dimensional totally real subspace of a space of type \,$(\Ieins,k')$\,; more explicitly: there exists \,$A \in \A$\,,
%a \,$2k'$-dimensional subspace \,$W \subset V(A)$\,, an orthogonal complex structure \,$j: W \to W$\, and a \,$k'$-dimensional totally real
%subspace \,$W'$\, of the unitary space \,$(W,j)$\, so that \,$U = \Menge{x+J(jx)}{x \in W'}$\,.
\end{itemize}
The various types of curvature-invariant spaces have the following properties:
\begin{center}
\begin{tabular}{|c|c|c|c|c|c|}
\hline
type of \,$\liem'$\, & $\dim_{\R} \liem'$ & \begin{minipage}{2cm} \begin{center} {\tiny \,$\liem'$\, complex or \\ totally real? \par} \end{center} \end{minipage} & $\rk(\liem')$ & \,$\liem'$\, maximal?\\
\hline
\hline
$(\Geo)$ & $1$ & totally real & $1$ & no \\ % , unless \,$m=1$\, \\
\hline
$(\Geins,k)$ & $2k$ & complex & $2$ & for \,$k = m-1 \geq 2$\, \\
$(\Gzwei,k_1,k_2)$ & \,$k_1+k_2$ & totally real & $2$ & for \,$k_1+k_2 = m \geq 3$\, \\
$(\Gdrei)$\, & $3$ & neither & $2$ & only for \,$m=2$\, \\
\hline
$(\Peins,k)$ & $k$ & totally real & $1$ & for \,$k = m$\, \\
$(\Pzwei)$ & $2$ & complex & $1$ & only for \,$m=2$\, \\
\hline
$(\Atyp)$ & $2$ & neither & $1$ & only for \,$m=3$\, \\
\hline
$(\Ieins,k)$ & $2k$ & complex & $1$ & for \,$2k = m \geq 4$ \\
$(\Izwei,k)$ & $k$ & totally real & $1$ & no \\
\hline
\end{tabular}
\end{center}
\end{Theorem}

\begin{Remark}
In the type specifications, the abbreviation ``Geo'' obviously stands for ``geodesic'', as the Lie triple systems of this type correspond to the
traces of geodesics in \,$Q$\,. The letters G, P, A and I stand for the words ``generic'' (because such spaces contain entire Cartan subalgebras of \,$Q$\,),
``principal'', ``arctan'' (because such spaces bear a relation to the angle \,$\arctan(\tfrac{1}{2})$\, as will be described in Section~\ref{SSe:cla:rk1})
and ``isotropic'', respectively.
\end{Remark}

\emph{Proof of Theorem~\ref{T:cla:cla}.} Via Proposition~\ref{P:Q:AR}(b) one verifies that the spaces mentioned in the theorem are indeed
curvature-invariant and therefore Lie triple systems, and one also sees easily that the information in the table on the dimension of the spaces
and on their being complex or totally real is correct. For the data on the rank of \,$\liem'$\,: Because the complex quadric has rank 2, we have 
\,$\rk(\liem') \in \{1,2\}$\, in any case, and \,$\rk(\liem') = 2$\, holds if and only if \,$\liem'$\, contains a Cartan subalgebra of \,$\liem$\,.
Via the explicit description of the Cartan subalgebras of \,$\liem$\, in Proposition~\ref{P:Q:roots}(a), one sees by this argument that 
\,$\liem'$\, is of rank 2 if it is of one of the types \,$(\Geins,k)$\,, \,$(\Gzwei,k_1,k_2)$\, or \,$(\Gdrei)$\,, and of rank 1 if it is of
any other type given in the theorem.

To prove the statements in the table on the maximality of Lie triple systems, we presume that the list of Lie triple systems given in the
theorem is complete; this fact will be proved in the remainder of the section. We then consider the various types individually:

\begin{itemize}
\item[$(\Geo,t)$]
If \,$\liem'$\, is of type \,$(\Geo,t)$\,, then \,$\liem'$\, is contained in a Cartan subalgebra, i.e. in a Lie triple system of type \,$(\Gzwei,1,1)$\,
and therefore cannot be maximal.
\item[$(\Geins,k)$]
This type exists only for \,$m \geq 3$\,.
If \,$\liem'$\, is of type \,$(\Geins,k)$\, with \,$k \leq m-2$\,, then \,$\liem'$\, is contained in a space of type \,$(\Geins,m-1)$\, and therefore cannot be maximal.
On the other hand, the spaces of type \,$(\Geins,m-1)$\, are of real codimension \,$2$\, in \,$\liem$\,. There exist no
Lie triple systems of \,$\liem$\, of real codimension \,$1$\, because of \,$m \geq 3$\,, and therefore the spaces of type \,$(\Geins,m-1)$\, are then maximal.
\item[$(\Gzwei,k_1,k_2)$]
If \,$\liem'$\, is of type \,$(\Gzwei,k_1,k_2)$\, with \,$k_1+k_2 < m$\,, then \,$\liem'$\, is contained in a space of type \,$(\Gzwei,k_1,m-k_1)$\, and is therefore
not maximal.
Moreover, any space \,$\liem'$\, of type \,$(\Gzwei,1,1)$\, is contained in a space of type \,$(\Gdrei)$\, and is therefore
not maximal in the case \,$m=2$\,. On the other hand, if \,$\liem'$\, is of type \,$(\Gzwei,k_1,k_2)$\, with \,$k_1+k_2 = m \geq 3$\,, then \,$\liem'$\, is maximal:
Assume to the contrary
that there exists a Lie triple system \,$\wh{\liem}'$\, of \,$\liem$\, with \,$\liem' \subsetneq \wh{\liem}' \subsetneq \liem$\,. Then we have
\,$\dim_{\R} \wh{\liem}' > \dim_{\R} \liem' = m$\,, and therefore \,$\wh{\liem}'$\, is of type \,$(\Geins,k)$\, for some \,$k$\, (see the table in the theorem) 
and hence complex.
Thus we have \,$\wh{\liem}' \supset \liem' \operp J\liem' = \liem$\,, which is a contradiction.
\item[$(\Gdrei)$]
For \,$m=2$\,, the spaces of type \,$(\Gdrei)$\, have real codimension \,$1$\, in \,$\liem$\, and are therefore maximal. On the other hand,
for \,$m \geq 3$\,, the space \,$\liem'$\, of type \,$(\Gdrei)$\, described in the theorem is contained in
the space \,$\C(x-Jy) \oplus \C(x+Jy) = \C x \oplus \C y$\, of type \,$(\Geins,2)$\,, and therefore cannot be maximal.
\item[$(\Peins,k)$]
If \,$\liem'$\, is of type \,$(\Peins,k)$\, with \,$k < m$\,, then \,$\liem'$\, is contained in a space of type \,$(\Peins,m)$\, and therefore cannot be maximal.
On the other hand, if \,$\liem'$\, is of type \,$(\Peins,m)$\,, then we have \,$\liem' = V(A)$\, for some \,$A \in \A$\,. An inspection of the table in the theorem
shows that there exists no Lie triple system \,$\wh{\liem}'$\, of \,$\liem$\, with \,$V(A) \subsetneq \wh{\liem}' \subsetneq \liem$\,.
\item[$(\Pzwei)$]
Let \,$\liem'$\, be a Lie triple system of type \,$(\Pzwei)$\,. In the case \,$m = 2$\,, \,$\liem'$\, is maximal: Assume to the contrary that there
exists a Lie triple system \,$\wh{\liem}'$\, of \,$\liem$\, with \,$\liem' \subsetneq \wh{\liem}' \subsetneq \liem$\,. Then \,$\wh{\liem}'$\, is of real dimension \,$3$\,
and therefore of type \,$(\Gdrei)$\,, so that there exists an orthonormal system \,$(x,y)$\, in some \,$V(A)$\,, \,$A \in \A$\, with
\,$\wh{\liem}' = \C(x-Jy) \operp \R(x+Jy)$\,.
\,$\liem'$\, is complex, and therefore we have \,$\liem' = \liem' \cap J\liem' \subset \wh{\liem}' \cap J\wh{\liem}' = \C(x-Jy)$\,, which is a contradiction because
\,$\C(x-Jy)$\, is an isotropic subspace of \,$\liem$\,, whereas \,$\liem'$\, is not.

On the other hand, in the case \,$m \geq 3$\,, \,$\liem'$\, is contained in a space of type \,$(\Geins,2)$\, and therefore cannot be maximal.
\item[$(\Atyp)$]
Let \,$\liem'$\, be a Lie triple system of type \,$(\Atyp)$\,, then we necessarily have \,$m \geq 3$\,.
Using the notation in the definition of this type in the theorem, we see that \,$\liem'$\, is contained in the Lie triple system
\,$\liem'' := \C x \oplus \C y \oplus \C z$\, of type \,$(\Geins,3)$\,; therefore \,$\liem'$\, cannot be maximal for \,$m \geq 4$\,. 

In the case \,$m = 3$\,, we again show
the maximality of \,$\liem'$\, by contradiction: Assume that \,$\wh{\liem}'$\, is a Lie triple system of \,$\liem$\, with
\,$\liem' \subsetneq \wh{\liem}' \subsetneq \liem$\,. Because \,$\liem'$\, contains vectors which are neither isotropic nor contained in some \,$V(A)$\,,
\,$\wh{\liem}'$\, is of one of the types \,$(\Geins,k)$\,, \,$(\Gzwei,k_1,k_2)$\,
and \,$(\Gdrei)$\,. If \,$\wh{\liem}'$\, is of type \,$(\Geins,k)$\,, then \,$\wh{\liem}'$\, is a \CQ-subspace of \,$\liem$\, and therefore contains \,$\liem''$\,;
because we have \,$\dim_{\C} \liem'' = 3 = \dim_{\C} \liem$\,, \,$\wh{\liem}' = \liem$\, follows, a contradiction. If \,$\wh{\liem}'$\, is of type \,$(\Gzwei,k_1,k_2)$\,,
then \,$\wh{\liem}'$\, is totally real in \,$\liem$\,, and hence \,$\liem'$\, is totally real, also a contradiction.
Finally, if \,$\wh{\liem}'$\, is of type \,$(\Gdrei)$\,, one obtains a contradiction to \,$\liem' \subset \wh{\liem}'$\,. 

\item[$(\Ieins,k)$]
Let \,$\liem'$\, be a Lie triple system of type \,$(\Ieins,k)$\,. Proposition~\ref{P:Q:CQ}(c)(i) shows that \,$\liem'$\, is properly contained in a \CQ-subspace
\,$\liem''$\, of \,$\liem$\, of complex dimension \,$2k$\,. 
In the case \,$2k<m$\, \,$\liem''$\, is  a Lie triple system of \,$\liem$\, of type \,$(\Geins,2k)$\,;
because we have \,$\wh{U} \supsetneq \liem'$\, it follows that \,$\liem'$\, is not maximal.
In the case \,$2k=m=2$\,, \,$\liem'$\, is contained in a space of type \,$(\Gdrei)$\, and therefore not maximal either.
In the case \,$2k=m\geq 4$\,, we once again prove the maximality of \,$\liem'$\, by contradiction:
Assume that \,$\wh{\liem}'$\, is a Lie triple system of \,$\liem$\, with \,$\liem' \subsetneq \wh{\liem}' \subsetneq \liem$\,. Then we have \,$\dim_{\R} \wh{\liem}'
> \dim_{\R} \liem' = 2k = m$\,, and therefore \,$\wh{\liem}'$\, is of type \,$(\Geins,k')$\, for some \,$k'$\,, and hence a \CQ-subspace. Thus we have
\,$\liem'' \subset \wh{\liem}'$\,; because of \,$\dim_{\C}(\liem'') = 2k = m$\,, we have \,$\liem'' = \liem$\, and therefore \,$\wh{\liem}' = \liem$\, follows,
a contradiction.
\item[$(\Izwei,k)$]
If \,$\liem'$\, is of type \,$(\Izwei,k)$\,, then \,$\liem'$\, is contained in the space \,$\liem' \operp J\liem'$\, of type \,$(\Ieins,k)$\, 
(see Proposition~\ref{P:Q:CQ}(c)(iii)) and therefore cannot be maximal.
\end{itemize}

It remains to prove that every Lie triple system in \,$\liem$\, is of one of the types given in the theorem, and this is the objective of
the remainder of the present section.

We make one preliminary observation:

\begin{Lemma}
\label{L:cla:complex}
Let \,$\liem'$\, be a Lie triple system of \,$\liem$\,, and suppose that there exist \,$A \in \A$\, and \,$X,Y \in V(A) \setminus \{0\}$\, with
\,$JX,Y \in \liem'$\, and \,$\g{X}{Y} \neq 0$\,. Then \,$\liem'$\, is a complex linear subspace of \,$\liem$\,. 
\end{Lemma}

\beweis
Let \,$R$\, be the curvature tensor of \,$Q$\, and let us fix \,$A \in \A$\,. Then we have by Proposition~\ref{P:Q:AR}(b) for any \,$Z \in \liem'$\,
\begin{align*}
R(JX,Y)Z & = \gC{Z}{Y}JX - \gC{Z}{JX}Y + 2\g{X}{Y}JZ + \gC{Y}{AZ}AJX - \gC{JX}{AZ}AY \\
& = \gC{Z}{Y}JX - \gC{Z}{JX}Y + 2\g{X}{Y}JZ + \gC{Z}{AY}AJX - \gC{Z}{AJX}AY \\
& = 2\g{X}{Y}JZ \; ; 
\end{align*}
in this calculation it has been used that the \,$\C$-bilinear form \,$\liem \times \liem \to \C,\;(v,w) \mapsto \gC{v}{Aw}$\, is symmetric
and that \,$AJX = -JX$\,, \,$AY = Y$\, holds. Because \,$\liem'$\, is a Lie triple system and thus curvature invariant, it follows that
\,$2\g{X}{Y}JZ \in \liem'$\, and hence \,$JZ \in \liem'$\, holds.
\beweisende

\subsection{The case of rank 2}
\label{SSe:cla:rk2}

We now start with the proof that the list of Lie triple systems given in Theorem~\ref{T:cla:cla} is in fact complete.

We let a Lie triple system \,$\{0\} \neq \liem' \subset \liem$\, be given. We have \,$\rk(\liem') \leq \rk(Q) = 2$\,
and therefore \,$\rk(\liem') \in \{1,2\}$\,. The two resulting cases \,$\rk(\liem') = 2$\, and \,$\rk(\liem') = 1$\, divide the proof of the
classification theorem into two main parts. We treat the case \,$\rk(\liem') = 2$\, in the present subsection, and the case \,$\rk(\liem') = 1$\, in the
next subsection.

Thus we now suppose \,$\rk(\liem') = 2$\,. We fix a Cartan subalgebra \,$\liea$\, of \,$\liem'$\,;
because of \,$\rk(\liem') = 2 = \rk(Q)$\,, \,$\liea$\, also is a Cartan subalgebra of \,$\liem$\,. In the sequel, we denote by \,$\Delta := \Delta(\liem,\liea)$\,
and \,$\Delta' := \Delta(\liem',\liea)$\, the root systems of \,$\liem$\, resp.~of \,$\liem'$\, with respect to \,$\liea$\,. In this relation
we use the notations introduced in Section~\ref{Se:roots}. Then we have by Proposition~\ref{P:cla:subroots:subroots-neu}(b)
\begin{equation}
\label{eq:cla:rk2:m'subsetm}
\Delta' \subset \Delta \qmq{and} \forall \alpha \in \Delta' \; : \; \liem_\alpha' = \liem_\alpha \cap \liem' \subset \liem_\alpha \; .
\end{equation}
Therefore \,$\Delta_+' := \Delta_+ \cap \Delta'$\, is a system of positive roots of \,$\Delta'$\,, where
\,$\Delta_+ := \{\lambda_1,\dotsc,\lambda_4\}$\, is  the system of positive roots of \,$\Delta$\, described in Proposition~\ref{P:Q:roots}(b).
Further, we have by Equation~\eqref{eq:roots:m'decomp}
\begin{equation}
\label{eq:cla:rk2:decomp}
\liem' = \liea \;\oplus\; \bigoplus_{\alpha \in \Delta_+'} \liem_\alpha' \; .
\end{equation}
% These relations are the pivotal point for the classification in the case \,$\rk(\liem') = 2$\,.

Moreover, the root system \,$\Delta'$\, is invariant under the Weyl group \,$W(\liem',\liea')$\, by Proposition~\ref{P:roots:weyl}(c), and this fact
imposes restrictions on the subsets of \,$\Delta_+ = \{\lambda_1,\dotsc,\lambda_4\}$\, which can occur as \,$\Delta_+'$\,. For example
\,$\Delta_+' = \{\lambda_1,\lambda_4\}$\, is impossible, because then \,$\Delta' =  \Delta_+' \dot{\cup} (-\Delta_+')$\, would not
be invariant under the reflection in the line orthogonal to \,$\lambda_1$\,.
%(For the calculation of the action of the Weyl group on the \,$\lambda_k$\,, note the relationship between its action on \,$\lambda_k$\,
%and on \,$\lambda_k^\sharp$\, given by Equation~\eqref{eq:Lts:roots:Weyl-lambda-lambdasharp} and the explicit description of the \,$\lambda_k^\sharp$\,
%in Theorem~\ref{T:iso:sym:roots}(b).)

By this consideration we see that \,$\Delta_+'$\, must be one of the following eight sets:
$$ \varnothing\;,\quad \{\lambda_1\}\;,\quad \{\lambda_2\}\;,\quad \{\lambda_3\}\;,\quad \{\lambda_4\}\;, \quad
\{\lambda_1,\lambda_2\}\;,\quad \{\lambda_3,\lambda_4\}\;,\quad \{\lambda_1,\lambda_2,\lambda_3,\lambda_4\} \; .
$$

We now inspect the eight cases of possible \,$\Delta_+'$\, individually to show that the corresponding Lie triple systems \,$\liem'$\,
are all of of one of the types \,$(\Geins,k)$\,, \,$(\Gzwei,k_1,k_2)$\, and \,$(\Gdrei)$\, as they are described in Theorem~\ref{T:cla:cla}.

For this purpose, we note that by Proposition~\ref{P:Q:roots}(a)
there exist \,$A \in \A$\, and an orthonormal system \,$(X,Y)$\, in \,$V(A)$\, so that \,$\liea = \R X \oplus \R JY$\, holds. Also, we put
\,$n_\alpha' := \dim (\liem_\alpha')$\, for \,$\alpha \in \Delta'$\,, and continually use the data on the root system \,$\Delta_+ = \{\lambda_1,\dotsc,\lambda_4\}$\,
and the root spaces \,$\liem_{\lambda_k}$\, given in Proposition~\ref{P:Q:roots}(b).

\emph{The case \,$\Delta_+' = \varnothing$\,.}
By Equation~\eqref{eq:cla:rk2:decomp} we have
\,$\liem' = \liea = \R X \oplus \R JY$\,, and therefore, \,$\liem'$\, is of type \,$(\Gzwei,1,1)$\,
with \,$W_1 := \R X$\,, \,$W_2 := \R Y$\,.

\emph{The case \,$\Delta_+' = \{\lambda_1\}$\,.}
By Equation~\eqref{eq:cla:rk2:decomp} we have
\,$\liem' = \liea \oplus \liem'_{\lambda_1}$\,; by \eqref{eq:cla:rk2:m'subsetm} we have
\,$\liem'_{\lambda_1} \subset \liem_{\lambda_1} = J((\R X \oplus \R Y)^{\perp,V(A)})$\,.
It follows that \,$\liem'$\, is of type \,$(\Gzwei,1,1+n'_{\lambda_1})$\, with \,$W_1 := \R X$\, and \,$W_2 := \R Y \operp J\liem'_{\lambda_1}$\,.

\emph{The case \,$\Delta_+' = \{\lambda_2\}$\,.}
Analogously as in the case \,$\Delta_+' = \{\lambda_1\}$\, we see that \,$\liem'$\, is of type \,$(\Gzwei,1+n'_{\lambda_2},1)$\, with
\,$W_1 := \R X \operp \liem'_{\lambda_2}$\, and \,$W_2 := \R Y$\,.

\emph{The case \,$\Delta_+' = \{\lambda_3\}$\,.}
By Equation~\eqref{eq:cla:rk2:decomp} we have \,$\liem' = \liea \oplus \liem_{\lambda_3}'$\,.
We have \,$\{0\} \neq \liem'_{\lambda_3} \subset \liem_{\lambda_3}$\,; because \,$\liem_{\lambda_3}$\, is 1-dimensional, therefrom already
\,$\liem'_{\lambda_3} = \liem_{\lambda_3} = \R(JX+Y)$\, follows. Thus we have
\begin{align*}
\liem' & = \liea \oplus \liem_{\lambda_3} = \R X \oplus \R JY \oplus \R(JX+Y) \\
& = \R(X+JY) \oplus \R(X-JY) \oplus \R(JX+Y) = \R(X+JY) \oplus \C(X-JY) \;,
\end{align*}
and therefore \,$\liem'$\, is of type \,$(\Gdrei)$\,. % with this choice of \,$X$\, and \,$Y$\,.

\emph{The case \,$\Delta_+' = \{\lambda_4\}$\,.}
Analogously as in the case \,$\Delta_+' = \{\lambda_3\}$\, we obtain \,$\liem' = \R(X-JY) \oplus \C(X+JY)$\,.
By replacing \,$Y$\, with \,$-Y$\,, we see that also in this case \,$\liem'$\, is of type \,$(\Gdrei)$\,.
%By replacing the conjugation \,$A \in \A$\, by \,$-A \in \A$\,, %one sees (because of \,$V(-A) = JV(A)$\,) that in this case \,$\liem'$\, also is of type \,$(\Gdrei)$\,.

\emph{The case \,$\Delta_+' = \{\lambda_1,\lambda_2\}$\,.}
By Equation~\eqref{eq:cla:rk2:decomp} we have
\begin{equation}
\label{eq:cla:rk2:2a-decomp}
\liem' = \liea \oplus \liem_{\lambda_1}' \oplus \liem_{\lambda_2}' = W_1 \oplus J(W_2)
\end{equation}
with \,$W_1 := \R X \oplus \liem_{\lambda_2}'$\, and \,$W_2 := \R Y \oplus J(\liem_{\lambda_1}')$\,.
Together with Equation~\eqref{eq:cla:rk2:m'subsetm},
the table in Proposition~\ref{P:Q:roots}(b) shows that \,$J(\liem_{\lambda_1}'), \liem_{\lambda_2}' \subset (\R X \oplus \R Y)^{\perp,V(A)} \subset V(A)$\, holds,
and therefore we have \,$W_1,W_2 \subset V(A)$\,. 

We now show \,$W_1 \perp W_2$\,: Let \,$u \in W_2$\, and \,$v \in W_1$\, be given, and assume that
\,$\g{u}{v} \neq 0$\, holds. We have \,$Ju,v \in \liem'$\, by Equation~\eqref{eq:cla:rk2:2a-decomp}, and therefore Lemma~\ref{L:cla:complex}
shows that \,$\liem'$\, is a complex-linear subspace of \,$\liem$\,. Because we have \,$X+JY \in \liea \subset \liem'$\,, it follows that we also have
\,$-Y+JX = J(X+JY) \in \liem'$\,.
Hence we have \,$\liem_{\lambda_4} = \R(JX-Y) \subset \liem'$\, (see Proposition~\ref{P:Q:roots}(b)) and therefore
\,$\liem_{\lambda_4}' = \liem_{\lambda_4} \cap \liem' = \liem_{\lambda_4}$\, (see Proposition~\ref{P:cla:subroots:subroots-neu}(b)),
whence \,$\lambda_4 \in \Delta_+'$\, follows.
But this is a contradiction to the hypothesis \,$\Delta_+' = \{\lambda_1,\lambda_2\}$\, defining the present case.

Therefore \,$\liem'$\, is of type \,$(\Gzwei,1+n_{\lambda_2}', 1+n_{\lambda_1}')$\, with the present choice of \,$W_1$\, and \,$W_2$\,.

\emph{The case \,$\Delta_+' = \{\lambda_3,\lambda_4\}$\,.}
For \,$k \in \{3,4\}$\, we have \,$\dim \liem_{\lambda_k} = 1$\,, and therefore
the same argument as in the treatment of the case \,$\Delta_+' = \{\lambda_3\}$\, shows that \,$\liem_{\lambda_k}' = \liem_{\lambda_k}$\, holds.
Thus we have by Equation~\eqref{eq:cla:rk2:decomp}
\begin{align*}
\liem' & = \liea \oplus \liem_{\lambda_3}' \oplus \liem_{\lambda_4}' = (\R X \oplus \R JY) \oplus \R(JX+Y) \oplus \R(JX-Y) \\
& = \R X \oplus \R JY \oplus \R JX \oplus \R Y = \C X \oplus \C Y \; .
\end{align*}
Thus we have \,$\liem' = W \oplus JW$\, with \,$W := \R X \operp \R Y \subset V(A)$\,. Therefore \,$\liem'$\, is a \,$2$-dimensional \CQ-subspace and hence of type
\,$(\Geins,2)$\,.

\emph{The case \,$\Delta_+' = \{\lambda_1,\lambda_2,\lambda_3,\lambda_4\}$\,.}
By Equation~\eqref{eq:cla:rk2:decomp} we have
\begin{equation}
\label{eq:cla:rk2:4-decomp}
\liem' = \liea \oplus \liem_{\lambda_1}' \oplus \liem_{\lambda_2}' \oplus \liem_{\lambda_3}' \oplus \liem_{\lambda_4}' \;,
\end{equation}
and by an analogous argument as for the case \,$\Delta_+' = \{\lambda_3,\lambda_4\}$\,, we see that
\begin{equation}
\label{eq:cla:rk2:4-a34}
\liem' \;\overset{\eqref{eq:cla:rk2:4-decomp}}{\supset}\; \liea \oplus \liem_{\lambda_3}' \oplus \liem_{\lambda_4}' = \C X \oplus \C Y
\end{equation}
holds. In particular we have \,$X,JX \in \liem'$\,, whence it follows by Lemma~\ref{L:cla:complex} that \,$\liem'$\, is a complex linear
subspace of \,$\liem$\,. Therefrom \,$\liem'_{\lambda_1} = J(\liem'_{\lambda_2})$\, follows, and thus we obtain from Equations~\eqref{eq:cla:rk2:4-decomp}
and \eqref{eq:cla:rk2:4-a34}:
$$ \liem' = \C X \oplus \C Y \oplus J(\liem'_{\lambda_2}) \oplus \liem'_{\lambda_2} = W \oplus JW $$
with \,$W := \R X \operp \R Y \operp \liem'_{\lambda_2} \subset V(A)$\,. Therefore \,$\liem'$\, is a \,$(2+n_{\lambda_2}')$-dimensional \CQ-subspace
and hence of type \,$(\Geins,2+n_{\lambda_2}')$\,.

This completes the classification of the rank \,$2$\, Lie triple systems in \,$\liem$\,.

\subsection{The case of rank 1}
\label{SSe:cla:rk1}

We first give a way to describe the position of vectors \,$v \in \liem$\, in the Cartan subalgebra in which they lie.
For this, let \,$v \in \liem \setminus \{0\}$\, be given. Then there exists a Cartan subalgebra \,$\liea \subset \liem$\, with \,$v \in \liea$\,
(see \cite{Loos:1969-2}, Theorem VI.1.2(b), p.~56), and by Proposition~\ref{P:Q:roots}(a) there exist \,$A \in \A$\, and \,$X,Y \in \Sph(V(A))$\, with 
\,$\liea = \R X \oplus \R JY$\,. Because the Weyl group of \,$\liem$\, acts transitively on the set of Weyl chambers in \,$\liea$\,
(\cite{Loos:1969-2}, Theorem VI.2.2, p.~67), there exists \,$g \in K$\, so that \,$\Ad(g)v$\, lies in the closed Weyl chamber of \,$\liea$\,
bounded by the root vectors \,$\lambda_2^\sharp$\, and \,$\lambda_4^\sharp$\, as defined in Proposition~\ref{P:Q:roots}(b):

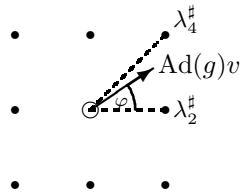
\begin{center}
\begin{minipage}{5cm}
\begin{center}
\strut \\[.3cm]
\setlength{\unitlength}{1cm}
\begin{picture}(2,2)
\put(1,1){\circle{0.2}}
\put(2,1){\circle*{0.1}}        % lambda_2
\put(2,0){\circle*{0.1}}        % lambda_3
\put(1,0){\circle*{0.1}}        % -lambda_1
\put(0,0){\circle*{0.1}}        % -lambda_4
\put(0,1){\circle*{0.1}}        % -lambda_2
\put(0,2){\circle*{0.1}}        % -lambda_3
\put(1,2){\circle*{0.1}}        % lambda_1
\put(2,2){\circle*{0.1}}        % lambda_4
%\put(2.1,-0.3){{\small $\lambda_3^\sharp$}}
\put(2.1,0.90){{\small $\lambda_2^\sharp$}}
\put(2.1,2.1){{\small $\lambda_4^\sharp$}}
% \put(0.6,2.15){{\small $\lambda_1^\sharp$}}
%\put(0.7,-0.4){{$-\lambda_1^\sharp$}}
%\put(-0.6,-0.3){{$-\lambda_4^\sharp$}}
%\put(-0.6,0.9){{$-\lambda_2^\sharp$}}
%\put(-0.6,2.1){{$-\lambda_3^\sharp$}}
\thicklines
\put(1,1){\vector(3,2){0.8333}}
\dashline[100]{0.1}(1,1)(2,1)
\dashline[100]{0.1}(1,1)(2,2)
%\dashline[50]{0.1}(1,1)(0,1)
%\dashline[50]{0.1}(1,1)(1,2)
%\dashline[50]{0.1}(1,1)(1,0)
\put(1,1){\arc(0.6,0){33.7}}
\put(1.32,1.05){{\scriptsize $\vi$}}
\put(1.9,1.5){$\Ad(g)v$}
\end{picture}
\strut \\[1cm]
\end{center}
\end{minipage}
\end{center}

In this situation, we call the angle \,$\vi \in [0,\tfrac{\pi}{4}]$\, between the vectors \,$\lambda_2^\sharp = \sqrt{2}\,X$\, and \,$\Ad(g)v$\, 
the \emph{characteristic angle} of \,$v$\, and denote it 
by \,$\vi(v)$\,. Because the action of the Weyl group on the set of Weyl chambers is simply transitive, this angle does not depend on the choice
of \,$g$\, (with the aforementioned property).

\begin{Prop}
\label{P:rk1:angle}
Suppose \,$v \in \liem \setminus \{0\}$\,. 
\begin{enumerate}
\item
The angle \,$\vi(v) \in [0,\tfrac{\pi}{4}]$\, is characterized by
$$ |\gC{v}{Av}| = \cos(2\vi(v)) \cdot \|v\|^2 \;, $$
where the number \,$|\gC{v}{Av}|$\, does not depend on the choice of \,$A \in \A$\,. 
Therefore \,$\vi(v)$\, is uniquely determined even in those cases where \,$v$\, is contained in more than one Cartan subalgebra.
\item
We have \,$\vi(cv) = \vi(v)$\, for any \,$c \in \R \setminus \{0\}$\,.
%\item
%The map \,$\vi: \liem \setminus \{0\},\; v \mapsto \vi(v)$\, is continuous.
\item
\,$\vi$\, is invariant under \CQ-automorphisms of \,$\liem$\,, in particular \,$\vi$\, is \,$\Ad(K)$-invariant.
\item
Let \,$\liea$\, be a Cartan subalgebra of \,$\liem$\, with \,$v \in \liea$\, and \,$v' \in \liea \setminus \{0\}$\, with \,$\g{v}{v'} = 0$\,. Then we have
\,$\vi(v) = \vi(v')$\,.
\end{enumerate}
\end{Prop}

\beweis
We first note that the map \,$f: \liem \to \R,\; v \mapsto |\gC{v}{Av}|$\, does not depend on the choice of \,$A \in \A$\, and is invariant under
\CQ-automorphisms, in particular (see Equation~\eqref{eq:iso:sym:Ad-tau} and Proposition~\ref{P:Q:Qaction}(b)) it is \,$\Ad(K)$-invariant.

Now let \,$v \in \liem \setminus \{0\}$\, be given; without loss of generality we may suppose \,$\|v\| = 1$\,. Let \,$\liea = \R X \oplus \R JY$\,
(\,$X,Y \in \Sph(V(A))$\,, \,$A \in \A$\,) be a Cartan subalgebra of \,$\liem$\, with \,$v \in \liea$\,, and let \,$g \in K$\, be so that
\,$\Ad(g)v$\, lies in the closed Weyl chamber of \,$\liea$\, bounded by \,$\lambda_2^\sharp$\, and \,$\lambda_4^\sharp$\,. Then we have
$$ \Ad(g)v = \cos(\vi(v))X + \sin(\vi(v))JY $$
and therefore
\begin{align*}
f(v) & \,=\, f(\Ad(g)v) = \bigr|\; \gC{\cos(\vi(v))X + \sin(\vi(v))JY}{\cos(\vi(v))X - \sin(\vi(v))JY} \;\bigr| \\
& \overset{(*)}{=} |\cos^2(\vi(v)) - \sin^2(\vi(v))| = |\cos(2\vi(v))| \;; 
\end{align*}
for the equals sign marked $(*)$ it should be noted that \,$\gC{X}{X} = \gC{Y}{Y} = 1$\, and \,$\gC{X}{Y} = 0$\, holds. This proves (a).
(b) is obvious and (c) follows from (a) via the fact that \,$f$\, is invariant under \CQ-automorphisms. 
%It is also a consequence of (a) that we
%have
%$$ \forall v \in \liem \setminus \{0\} \; : \; \vi(v) = \tfrac{1}{2} \arccos \left( \frac{|\gC{v}{Av}|}{\|v\|^2} \right) \; ; $$
%from this representation of \,$\vi$\,, (c) follows.

In the situation of (d), %by Proposition~\ref{P:Q:roots}(a) there exist
% \,$A \in \A$\, and an orthonormal system \,$(x,y)$\, in \,$V(A)$\, so that \,$\liea = \R x \oplus \R Jy$\, holds.
% Then there exists \,$t \in \R$\, so that \,$v = \cos(t)x + \sin(t)Jy$\, holds, and with the same \,$t$\,, we have \,$v' = \pm(-\sin(t)x + \cos(t)Jy)$\,. 
we consider the Weyl transformation \,$B: \liea \to \liea$\, obtained by first reflecting in the line perpendicular to \,$\lambda_1^\sharp$\, and then
reflecting in the line perpendicular to \,$\lambda_4^\sharp$\,. Then \,$B$\, is a rotation in the plane \,$\liea$\, by the angle \,$\tfrac{\pi}{2}$\,,
and therefore there exists \,$c \in \R \setminus \{0\}$\, so that \,$v' = c\cdot Bv$\, holds. We now have
$$ \vi(v') = \vi(c\cdot Bv) \overset{(b)}{=} \vi(Bv) \overset{(*)}{=} \vi(v) \; ; $$
for the equals sign marked $(*)$, see Proposition~\ref{P:roots:weyl}(b) and part (c) of the present proposition.
\beweisende

\begin{Remark}
Via a somewhat different approach, such a characteristic angle has already been introduced in \cite{Reckziegel:quadrik-1995}. 
In Section~6 of \cite{Reckziegel:quadrik-1995} it is shown that the orbits of the isotropy action of \,$\SO(2) \times \SO(m)$\, on the unit sphere \,$\Sph(T_pQ)$\,
(see Proposition~\ref{P:Q:Qaction}(a),(b))
are exactly the sets \,$M_t := \Menge{v \in \Sph(T_pQ)}{\vi(v) = t}$\, with \,$t \in [0,\tfrac{\pi}{4}]$\,. The orbits of the isotropy
actions of rank-2-symmetric spaces on the unit sphere have been studied extensively by \textsc{Takagi} and \textsc{Takahashi} in 
\cite{Takagi/Takahashi:homhyp-1972}. There it is shown that \,$(M_t)_{0<t<\pi/4}$\, is a
family of isoparametric hypersurfaces in \,$\Sph(T_pQ)$\,. It has \,$g=4$\, principal curvatures for \,$m \geq 3$\,; for \,$m=2$\, the number of
principal curvatures is reduced to \,$g=2$\,. The focal sets of \,$(M_t)_{0<t<\pi/4}$\, are \,$M_0$\, and \,$M_{\pi/4}$\,. 
\end{Remark}

\begin{Prop}
\label{P:rk1:vi0}
Let \,$\liem'$\, be a Lie triple system in \,$\liem$\, of rank \,$1$\,. 
Then all elements of \,$\liem' \setminus \{0\}$\, have one and the same characteristic angle \,$\vi_0 \in [0,\tfrac{\pi}{4}]$\,. If \,$\dim(\liem') \geq 2$\, holds,
then we have \,$\vi_0 \in \{0,\arctan(\tfrac{1}{2}),\tfrac{\pi}{4}\}$\,, and in the case \,$\vi_0 = \arctan(\tfrac{1}{2})$\,, \,$\liem'$\, does not
have any elementary roots (see Definition~\ref{D:cla:subroots:Elemcomp}). 
\end{Prop}

\beweis
We consider the Lie subalgebra \,$\liek' := [\liem',\liem']$\, of \,$\liek$\, and the connected Lie subgroup \,$K'$\, of \,$K$\, with 
Lie algebra \,$\liek'$\,. We let \,$Z_1,Z_2 \in \liem' \setminus \{0\}$\, be given. Then \,$\R Z_1$\, and \,$\R Z_2$\, are Cartan subalgebras of \,$\liem'$\,,
and therefore there exists \,$g \in K'$\, with \,$\Ad(g)\R Z_1 = \R Z_2$\, by \cite{Loos:1969-2}, Theorem~VI.1.2(c) (note that \,$\liem'$\, has to be 
of compact type in this situation).
Hence there exists \,$c \in \R \setminus \{0\}$\,
with \,$\Ad(g)Z_1 = cZ_2$\,, and therefore we have by Proposition~\ref{P:rk1:angle}(c),(b)
$$ \vi(Z_1) = \vi(\Ad(g)Z_1) = \vi(cZ_2) = \vi(Z_2) \; . $$
This proves the first statement of the proposition. 

We now suppose that \,$\dim(\liem') \geq 2$\, holds and let any \,$Z \in \liem' \setminus \{0\}$\, be given. \,$\liea' := \R Z$\, is a Cartan subalgebra of
\,$\liem'$\,; \,$\liea'$\, is contained in a Cartan subalgebra \,$\liea$\, of \,$\liem$\, by \cite{Loos:1969-2}, Theorem~VI.1.2(b). 
We now consider the root systems \,$\Delta' := \Delta(\liem',\liea')$\, and \,$\Delta := \Delta(\liem,\liea)$\, of \,$\liem'$\, resp.~\,$\liem$\,
with respect to \,$\liea'$\, resp.~\,$\liea$\,.
Because of \,$\dim(\liem') > \rk(\liem')$\,, we have \,$\Delta' \neq \varnothing$\,. We consider any root \,$\alpha \in \Delta'$\,. 
By Equation~\eqref{eq:cla:subroots:subroots-neu:toshow-Delta} in Proposition~\ref{P:cla:subroots:subroots-neu}(a), \,$\alpha$\, has to be either
elementary or composite in the sense of Definition~\ref{D:cla:subroots:Elemcomp}.

If \,$\alpha$\, is elementary and \,$\lambda \in \Delta$\, is the unique root with \,$\lambda|\liea' = \alpha$\,, then \,$Z$\, is co-linear to \,$\lambda^\sharp$\,
by Proposition~\ref{P:cla:subroots:Comp}(a), and therefore we have \,$\vi_0 = \vi(Z) = \vi(\lambda^\sharp)$\,. Via Proposition~\ref{P:rk1:angle}(a)
and the explicit representations of the root vectors \,$\lambda_k^\sharp$\, in Proposition~\ref{P:Q:roots}(b)  one easily calculates
\,$\vi(\lambda_1^\sharp) = \vi(\lambda_2^\sharp) = 0$\, and \,$\vi(\lambda_3^\sharp) = \vi(\lambda_4^\sharp) = \tfrac{\pi}{4}$\,.
Because we have \,$\lambda \in \Delta = \{\pm\lambda_1^\sharp,\dotsc,\pm\lambda_4^\sharp\}$\,, we thus see that 
\,$\vi_0 \in \{0,\tfrac{\pi}{4}\}$\, holds in the present case.

Now let us suppose that \,$\alpha$\, is composite, and let \,$\lambda,\mu \in \Delta$\, be two roots with \,$\lambda|\liea' = \alpha = \mu|\liea'$\, and
\,$\mu \neq \lambda$\,;
we also have \,$\mu \neq -\lambda$\, (because otherwise we would have \,$\alpha = 0 \not\in \Delta'$\,). By Proposition~\ref{P:cla:subroots:Comp}(b)
we have
$$ \g{Z}{\mu^\sharp - \lambda^\sharp} = 0 \;, $$
whence
\begin{equation}
\label{eq:cla:rk1:geo:vi}
\vi_0 = \vi(Z) = \vi(\mu^\sharp - \lambda^\sharp)
\end{equation}
follows by Proposition~\ref{P:rk1:angle}(d).

We now denote by \,$W$\, the Weyl group of \,$\Delta$\,, see Definition~\ref{D:roots:weyl}. Then we have for any \,$B \in W$\,
\begin{equation}
\label{eq:cla:rk1:geo:g-neq}
B(\lambda^\sharp) \neq \pm \, B(\mu^\sharp)
\end{equation}
and by Equation~\eqref{eq:cla:rk1:geo:vi} and Proposition~\ref{P:rk1:angle}(c)
\begin{equation}
\label{eq:cla:rk1:geo:g-vi}
\vi(Z) = \vi(\,B(\mu^\sharp) - B(\lambda^\sharp)\,) \; .
\end{equation}

\,$W$\, acts transitively on the complementary subsets \,$\{\pm \lambda_1^\sharp, \pm\lambda_2^\sharp\}$\, and \,$\{\pm \lambda_3^\sharp, \pm\lambda_4^\sharp\}$\,
of \,$\Delta$\,. Therefore there exists \,$B_1 \in W$\, so that \,$B_1(\lambda^\sharp) \in \{\lambda_1^\sharp,\lambda_3^\sharp\}$\, holds.

Let us first consider the case \,$B_1(\lambda^\sharp) = \lambda_1^\sharp$\,. Then we have \,$B_1(\mu^\sharp) \in
\{\pm \lambda_2^\sharp, \pm \lambda_3^\sharp, \pm \lambda_4^\sharp\}$\, by \eqref{eq:cla:rk1:geo:g-neq}. In the case
\,$B_1(\mu^\sharp) \in \{- \lambda_2^\sharp, - \lambda_3^\sharp, - \lambda_4^\sharp\}$\, we let \,$B_2 \in W$\, be the
reflection in \,$(\lambda_2^\sharp)^{\perp,\liea}$\, and have
$$ B_2(\lambda_1^\sharp) = \lambda_1^\sharp\,,\quad B_2(-\lambda_2^\sharp) = \lambda_2^\sharp\,,\quad
B_2(-\lambda_3^\sharp) = \lambda_4^\sharp \qmq{and} B_2(-\lambda_4^\sharp) = \lambda_3^\sharp \;; $$
otherwise we put \,$B_2 := \id_{\liea} \in W$\,. Thus with \,$B := B_2 \circ B_1 \in W$\, we have
\begin{equation}
\label{eq:cla:rk1:geo:g-l1-glm}
B(\lambda^\sharp) = \lambda_1^\sharp \qmq{and} B(\mu^\sharp) \in \{\lambda_2^\sharp, \lambda_3^\sharp, \lambda_4^\sharp\} \; .
\end{equation}
A calculation using Proposition~\ref{P:rk1:angle}(a) and the explicit presentation of \,$\lambda_k^\sharp$\, in Proposition~\ref{P:Q:roots}(b) shows
$$ \vi(\lambda_2^\sharp-\lambda_1^\sharp) = \tfrac{\pi}{4} \,,\quad \vi(\lambda_3^\sharp-\lambda_1^\sharp) = \arctan(\tfrac{1}{2}) \qmq{and}
\vi(\lambda_4^\sharp-\lambda_1^\sharp) = 0 \; . $$
In conjunction with Equations~\eqref{eq:cla:rk1:geo:g-vi} and \eqref{eq:cla:rk1:geo:g-l1-glm} these equations show that
(under the case hypotheses that \,$\alpha$\, is composite and \,$B_1(\lambda^\sharp) = \lambda_1^\sharp$\, holds) we have
\,$\vi(Z) \in \{0,\arctan(\tfrac{1}{2}),\tfrac{\pi}{4}\}$\,.

We now turn to the case \,$B_1(\lambda^\sharp) = \lambda_3^\sharp$\,. By an analogous argument as in the case \,$B_1(\lambda^\sharp) = \lambda_1^\sharp$\,
we see that there exists \,$B \in W$\, so that
\begin{equation}
\label{eq:cla:rk1:geo:g-l3-glm}
B(\lambda^\sharp) = \lambda_3^\sharp \qmq{and} B(\mu^\sharp) \in \{\lambda_1^\sharp, \lambda_2^\sharp, \lambda_4^\sharp\}
\end{equation}
holds, and a further calculation shows
$$ \vi(\lambda_1^\sharp-\lambda_3^\sharp) = \arctan(\tfrac{1}{2}) \,,\quad \vi(\lambda_2^\sharp-\lambda_3^\sharp) = 0 \qmq{and}
\vi(\lambda_4^\sharp-\lambda_3^\sharp) = 0 \; . $$
These facts together with Equations~\eqref{eq:cla:rk1:geo:g-vi} and \eqref{eq:cla:rk1:geo:g-l3-glm} show that in the present case we
have \,$\vi(Z) \in \{0,\arctan(\tfrac{1}{2})\}$\,.

This shows that in any case with \,$\dim(\liem')\geq 2$\,, \,$\vi_0 \in \{0,\arctan(\tfrac{1}{2}),\tfrac{\pi}{4}\}$\, holds. Moreover, we saw that if \,$\liem'$\, has an elementary root,
then in fact \,$\vi_0 \in \{0,\tfrac{\pi}{4}\}$\, holds.
\beweisende

We now classify the Lie triple systems \,$\liem' \subset \liem$\, of rank \,$1$\,. If \,$\dim(\liem') = 1$\, holds, then \,$\liem'$\, is of type
\,$(\Geo)$\,. Thus we may now suppose \,$\dim(\liem') \geq 2$\,. We fix \,$Z \in \Sph(\liem')$\,, then \,$\liea' := \R Z$\, is a Cartan subalgebra
of \,$\liem'$\,. As in the proof of Proposition~\ref{P:rk1:vi0} we choose a Cartan subalgebra \,$\liea$\, of \,$\liem$\, so that \,$\liea' = \liea \cap \liem'$\,
holds. By Proposition~\ref{P:Q:roots}(a) there exist \,$A \in \A$\, and an orthonormal system \,$(X,Y)$\, in \,$V(A)$\, so that \,$\liea = \R X \oplus \R JY$\,
holds, and by changing the sign of \,$A$\,, \,$X$\, and \,$Y$\, as necessary, we can ensure that \,$Z$\, lies in the closed Weyl chamber of \,$\liea$\, 
bounded by \,$X$\, and \,$X+JY$\,. Then we have
\begin{equation}
\label{eq:cla:rk1:ZXY}
Z = \cos(\vi_0)X + \sin(\vi_0)JY \; 
\end{equation}
with the constant \,$\vi_0 \in \{0,\arctan(\tfrac{1}{2}),\tfrac{\pi}{4}\}$\, from Proposition~\ref{P:rk1:vi0}.
Moreover, we consider the root systems \,$\Delta' := \Delta(\liem',\liea')$\, and \,$\Delta := \Delta(\liem,\liea)$\, of \,$\liem'$\, resp.~\,$\liem$\,
with respect to \,$\liea'$\, resp.~\,$\liea$\, and fix a system of positive roots \,$\Delta_+'$\, in \,$\Delta'$\,.
Then we have by Proposition~\ref{P:cla:subroots:subroots-neu}(a)
\begin{equation}
\label{eq:cla:rk1:Delta'}
\Delta' \subset \Menge{\lambda|\liea'}{\lambda \in \Delta,\,\lambda(Z)\neq 0}
\end{equation}
and
\begin{equation}
\label{eq:cla:rk1:decomp-m'}
\liem' = \R Z \;\oplus\; \bigoplus_{\alpha \in \Delta_+'} \liem_\alpha'
\end{equation}
with
\begin{equation}
\label{eq:cla:rk1:malpha'}
\forall \alpha \in \Delta_+' \; : \; \liem_\alpha' \;=\; \left( \bigoplus_{\substack{\lambda\in \Delta \\ \lambda(Z) = \alpha(Z)}} \liem_\lambda \right) 
\,\cap\, \liem'\; ;
\end{equation}
moreover we have because of \,$\dim(\liem') > \rk(\liem')$\,
\begin{equation}
\label{eq:cla:rk1:Delta'-not-empty}
\Delta' \neq \varnothing \; ,
\end{equation}
and Proposition~\ref{P:rk1:vi0} shows that on \,$\liem' \setminus \{0\}$\,, the characteristic angle function \,$\vi$\, is 
equal to some constant \,$\vi_0 \in \{0,\arctan(\tfrac{1}{2}),\tfrac{\pi}{4}\}$\,.
To complete the classification, we now treat the three possible values for \,$\vi_0$\, individually.

\textbf{The case \,$\boldsymbol{\vi_0 = 0}$\,.}
Then we have \,$Z = X$\, by Equation~\eqref{eq:cla:rk1:ZXY}.
%By Theorem~\ref{T:iso:sym:roots}(a) there exists \,$A \in \A$\, and an orthonormal system \,$(X,Y)$\, in \,$V(A)$\, so that \,$\liea = \R X \oplus \R JY$\,
%holds; moreover, because of \,$\vi(Z) = 0$\, these data can be chosen such that \,$Z = X$\, holds.
By Proposition~\ref{P:Q:roots}(b) we have
$$ \lambda_1(X) = 0 \qmq{and} \lambda_2(X) = \lambda_3(X) = \lambda_4(X) = \sqrt{2} \; ; $$
therefrom we conclude by \eqref{eq:cla:rk1:Delta'} and \eqref{eq:cla:rk1:Delta'-not-empty}
$$ \Delta' = \{\pm \alpha\} \qmq{with} \alpha(tZ) = \sqrt{2}\cdot t \text{ for \,$t \in \R$\,} $$
and by \eqref{eq:cla:rk1:decomp-m'} and \eqref{eq:cla:rk1:malpha'}
\begin{equation}
\label{eq:cla:rk1:0:decomp}
\liem' = \R X \oplus \liem_{\alpha}' \qmq{with} \{0\} \neq \liem_{\alpha}' \subset \liem_{\lambda_2} \oplus \liem_{\lambda_3} \oplus \liem_{\lambda_4} \; .
\end{equation}

Immediately, we will show that
\begin{equation}
\label{eq:cla:rk1:0:eitheror}
\text{either} \quad \liem'_{\alpha} \subset (\R X)^{\perp,V(A)} \qmq{or} \liem'_{\alpha} = \R\cdot JX
\end{equation}
holds. Then we conclude: In the case \,$\liem'_{\alpha} \subset (\R X)^{\perp,V(A)}$\,
we have \,$\liem' = \liea' \oplus \liem'_{\alpha} \subset V(A)$\,, therefore \,$\liem'$\, is of type \,$(\Peins,1+\dim \liem'_{\alpha})$\,.
On the other hand, in the case \,$\liem'_{\alpha} = \R \cdot JX$\, we have \,$\liem' = \liea' \oplus \liem'_{\alpha} = \C X$\,,
therefore \,$\liem'$\, is of type \,$(\Pzwei)$\,.

We now prove \eqref{eq:cla:rk1:0:eitheror}: Let \,$H \in \liem'_\alpha$\, be given. Then we have by \eqref{eq:cla:rk1:0:decomp} and
Proposition~\ref{P:Q:roots}(b)
$$ H \in  \liem_{\lambda_2} \oplus \liem_{\lambda_3}
                \oplus \liem_{\lambda_4}
        = \R\cdot JX \oplus (\R X)^{\perp, V(A)}
$$
and therefore there exist
\,$t \in \R$\, and \,$X' \in V(A)$\, with \,$X' \perp X$\, so that \,$H = t\cdot JX + X'$\, holds.
Via Proposition~\ref{P:Q:AR}(b) we calculate 
\begin{equation}
\label{eq:cla:rk1:0:tildew}
\widetilde{H} := \tfrac{1}{2} R(X,H)H  = (\|X'\|^2+t^2)\cdot X - 2t\cdot JX' \; .
\end{equation}
Because \,$\liem'$\, is curvature-invariant, we have \,$\widetilde{H} \in \liem'$\,. As \,$\liem'$\, is orthogonal to \,$\R JY \oplus \liem_{\lambda_1}
= (\R JX)^{\perp,JV(A)}$\, by Equation~\eqref{eq:cla:rk1:0:decomp} and hence in particular to \,$JX'$\,, we therefore have
$$ 0 = \g{\wt{H}}{JX'} = (-2t)\cdot \g{JX'}{JX'} = (-2t) \cdot \|X'\|^2 \; . $$
%we read off Equation~\eqref{eq:cla:rk1:0:tildew}
%that \,$\widetilde{H} \in \liea \oplus \liem_{\lambda_1}$\, holds (see Theorem~\ref{T:iso:sym:roots}(b)), and therefore we have
%\,$\widetilde{H} \in (\liea \oplus \liem_{\lambda_1}) \cap \liem' = \R X$\, by \eqref{eq:cla:rk1:0:decomp}.
%Equation~\eqref{eq:cla:rk1:0:tildew} thus implies \,$2t\cdot JX' = 0$\,.
Therefore we have either \,$t = 0$\,, implying \,$H = X' \in (\R X)^{\perp,V(A)}$\,; or else \,$\|X'\| = 0$\,, implying \,$H = t\cdot JX \in \R JX$\,. Thus, we have shown
$$ \liem'_{\alpha} \subset (\R X)^{\perp,V(A)} \;\cup\; \R\cdot JX \; . $$
Because \,$\liem'_{\alpha}$\, is a linear space, we in fact have
$$ \text{either} \quad \liem'_{\alpha} \subset (\R X)^{\perp,V(A)} \qmq{or} \liem'_{\alpha} \subset \R\cdot JX \;; $$
if the second case holds, then we actually have \,$\liem'_{\alpha} = \R\cdot JX$\, because of \,$\liem'_{\alpha} \neq \{0\}$\,.
Thus \eqref{eq:cla:rk1:0:eitheror} is shown.

\textbf{The case \,$\boldsymbol{\vi_0 = \arctan(\frac12)}$\,.}
By Equation~\eqref{eq:cla:rk1:ZXY} we then have
%As in the case \,$\vi_0 = 0$\,,
%by Theorem~\ref{T:iso:sym:roots}(a) there exists \,$A \in \A$\, and an orthonormal system \,$(X,Y)$\, in \,$V(A)$\, so that \,$\liea = \R X \oplus \R JY$\,
%holds; moreover, because of \,$\vi(Z) = \arctan(\frac12)$\, these data can be chosen such that
\begin{equation}
\label{eq:cla:rk1:atan12:v}
Z = \tfrac{2}{\sqrt{5}}X + \tfrac{1}{\sqrt{5}}JY \; ,
\end{equation}
and from Proposition~\ref{P:Q:roots}(b) we thus obtain
\begin{equation}
\label{eq:cla:rk1:atan12:lambdak}
\lambda_1(Z) = \tfrac{\sqrt{2}}{\sqrt{5}}\;,\quad \lambda_2(Z) = 2\,\tfrac{\sqrt{2}}{\sqrt{5}}\;,\quad \lambda_3(Z) = \tfrac{\sqrt{2}}{\sqrt{5}}
\qmq{and} \lambda_4(Z) = 3\,\tfrac{\sqrt{2}}{\sqrt{5}} \; .
\end{equation}
Because of \,$\vi_0 = \arctan(\tfrac{1}{2})$\, Proposition~\ref{P:rk1:vi0} shows that there do not exist any elementary roots in \,$\Delta'$\,;
therefore we conclude from Equations~\eqref{eq:cla:rk1:atan12:lambdak} by \eqref{eq:cla:rk1:Delta'} and \eqref{eq:cla:rk1:Delta'-not-empty}
$$ \Delta' = \{\pm \alpha\} \qmq{with} \alpha(tZ) = \tfrac{\sqrt{2}}{\sqrt{5}} \cdot t \text{ for \,$t \in \R$\,} $$
and by \eqref{eq:cla:rk1:decomp-m'} and \eqref{eq:cla:rk1:malpha'}
\begin{equation}
\label{eq:cla:rk1:atan12:decomp}
\liem' = \R Z \oplus \liem_{\alpha}' \qmq{with} \{0\} \neq \liem_{\alpha}' \subset \liem_{\lambda_1} \oplus \liem_{\lambda_3} \; .
\end{equation}

We now show
\begin{equation}
\label{eq:cla:rk1:atan12:m13}
\forall\, H \in \Sph(\liem'_{\alpha}) \; \exists\, U \in \Sph(V(A)) \; : \;
\bigr(\; H = \pm \tfrac{1}{\sqrt{5}}(Y + JX + \sqrt{3}JU) \qmq{and} U \perp X,Y \;\bigr) \; .
\end{equation}
Let \,$H \in \Sph(\liem'_{\alpha})$\, be given. We have \,$\vi(H) = \vi_0 = \arctan(\tfrac{1}{2})$\, and therefore by Proposition~\ref{P:rk1:angle}(a)
\begin{equation}
\label{eq:cla:rk1:atan12:wvi}
|\gC{H}{A(H)}| = \underbrace{\cos(2\,\arctan(\tfrac{1}{2}))}_{= \tfrac{3}{5}} \; . 
\end{equation}
By \eqref{eq:cla:rk1:atan12:decomp} and Proposition~\ref{P:Q:roots}(b) we have
$$ H \in \liem_{\lambda_1} \oplus \liem_{\lambda_3} = J(\R X \oplus \R Y)^{\perp,V(A)} \oplus \R(JX+Y) \; . $$
Consequently there exist \,$U' \in V(A)$\, with \,$U' \perp X,Y$\, and \,$t \in \R$\, so that
\begin{equation}
\label{eq:cla:rk1:atan12:w}
H = JU' + t\cdot (JX+Y) = tY + J(U' + tX)
\end{equation}
and therefore also
\begin{equation}
\label{eq:cla:rk1:atan12:Aw}
A(H) = tY - J(U' + tX)
\end{equation}
holds. Via Equations~\eqref{eq:cla:rk1:atan12:w} and \eqref{eq:cla:rk1:atan12:Aw} we obtain
$$ |\gC{H}{A(H)}| = \|U'\|^2 \;, $$
and by plugging the latter equation into \eqref{eq:cla:rk1:atan12:wvi}, we derive
\begin{equation*}
%\label{eq:cla:rk1:atan12:U'norm}
\|U'\|^2 = \tfrac{3}{5} \; . 
\end{equation*}
We now obtain from  Equation~\eqref{eq:cla:rk1:atan12:w} and the preceding equation 
$$ 1 = \|H\|^2 = 2t^2 + \|U'\|^2 = 2t^2 + \tfrac{3}{5} $$
Thus we have shown that
$$ t = \eps\, \tfrac{1}{\sqrt{5}} \qmq{and} \|U'\| = \sqrt{\tfrac{3}{5}} $$
holds with suitable \,$\eps \in \{\pm 1\}$\,. 
Consequently, we have \,$U := \eps\,\sqrt{{5}/{3}}\cdot U' \in \Sph(V(A))$\,. Equation~\eqref{eq:cla:rk1:atan12:w} shows that we have
\,$H = \eps\,\tfrac{1}{\sqrt{5}}(Y+JX+\sqrt{3}\,JU)$\,, and therefore \eqref{eq:cla:rk1:atan12:m13} is satisfied with this choice of \,$U$\,.

Next we prove \,$\dim \liem'_{\alpha} = 1$\,: Let \,$H_1, H_2 \in \Sph(\liem'_{\alpha})$\, be given; we will
show \,$H_2 = \pm H_1$\,. By \eqref{eq:cla:rk1:atan12:m13}, there exist \,$\eps_1, \eps_2 \in \{\pm 1\}$\, and
\,$U_1,U_2 \in \Sph(V(A))$\, so that
$$ H_k = \tfrac{\eps_k}{\sqrt{5}} \cdot (Y+JX+\sqrt{3}JU_k) $$
holds for \,$k \in \{1,2\}$\,. Under the assumption \,$H_2 \neq \pm H_1$\, we could suppose without loss of generality that \,$\eps_1 = \eps_2 = 1$\, holds,
and then \,$H_1-H_2 = \sqrt{3/5}\cdot J(U_1-U_2) \in \liem'_{\alpha} \subset \liem'$\, would be a non-zero vector; it is contained in \,$JV(A)$\, and
therefore has the property \,$\vi(H_1-H_2) = 0$\,, in contradiction to \,$\forall H \in \liem'\setminus\{0\}: \vi(H) = \arctan(\tfrac12)$\,.

Thus \,$\liem'_{\alpha}$\, is 1-dimensional, and therefore we have \,$\liem' = \liea' \oplus \liem'_{\alpha} = \R Z \oplus \R H$\, with any \,$H \in \Sph(\liem'_{\alpha})$\,.
Equations~\eqref{eq:cla:rk1:atan12:v} and \eqref{eq:cla:rk1:atan12:m13} now show that \,$\liem'$\, is a space of type \,$(\Atyp)$\,.

\textbf{The case \,$\boldsymbol{\vi_0 = \tfrac\pi4}$\,.}
In this case we have for every \,$Z \in \liem'\setminus \{0\}$\,: \,$\cos(2\,\vi(Z))=\cos(\tfrac{\pi}{2}) = 0$\, and therefore by
Proposition~\ref{P:rk1:angle}(a) \,$\gC{Z}{A(Z)} = 0$\,. This shows that \,$\liem'$\, is isotropic (see Definition~\ref{D:Q:CQ}(e)). Therefore
the ``complex closure'' \,$\wh{\liem}' := \liem' + J\liem' \subset \liem$\, of \,$\liem'$\, is also isotropic by Proposition~\ref{P:Q:CQ}(c)(iii),
and therefore a curvature-invariant subspace of \,$\liem$\, of type \,$(\Ieins,k)$\, with \,$k := \dim_{\C} \wh{\liem}'$\,.
Hence the quadratic form corresponding to the \,$\C$-bilinear form \,$\beta: \wh{\liem}' \times \wh{\liem}' \to \C,\; (Z_1,Z_2) \mapsto \gC{Z_1}{A(Z_2)}$\,
vanishes, and therefore we have \,$\beta=0$\,. From this fact and Proposition~\ref{P:Q:AR}(b) we see that for \,$Z_1,Z_2,Z_3 \in \wh{\liem}'$\,, the curvature
tensor \,$R$\, of \,$Q$\, is given by
\begin{align*}
R(Z_1,Z_2)Z_3 & = \gC{Z_3}{Z_2}Z_1 - \gC{Z_3}{Z_1}Z_2 - 2\cdot\g{JZ_1}{Z_2}JZ_3 \\
& = \g{Z_2}{Z_3}Z_1 - \g{Z_1}{Z_3}Z_2 + \g{JZ_2}{Z_3}JZ_1 - \g{JZ_1}{Z_3}JZ_2 - 2\cdot \g{JZ_1}{Z_2}JZ_3 \; . 
\end{align*}
Therefore the restriction of the curvature tensor of \,$Q$\, to \,$\wh{\liem}'$\, is 
the curvature tensor of a complex projective space of constant holomorphic sectional curvature~\,$4$\,.

If \,$\liem'$\, is a complex subspace of \,$\liem$\,, we have \,$\liem' = \widehat{\liem}'$\,; therefore \,$\liem'$\, then is of type \,$(\Ieins,k)$\,.
Otherwise, \,$\liem'$\, is a Lie triple system of \,$\widehat{\liem}'$\,; by the well-known classification of totally geodesic submanifolds
in a complex projective space, it follows that \,$\liem'$\, is a totally real subspace of \,$\widehat{\liem}'$\,, and therefore a \,$k$-dimensional
totally real, isotropic subspace of \,$\liem$\,. Consequently, \,$\liem'$\, is then of type \,$(\Izwei,k)$\,.

This completes the proof of Theorem~\ref{T:cla:cla}.
\hfill
$\Box$

\section{Totally geodesic embeddings into the complex quadric}
\label{Se:tgsub}

We now wish to find out which (connected, complete) totally geodesic submanifolds \,$M$\, of \,$Q$\, correspond to the Lie triple systems \,$\liem'$\, classified
in Theorem~\ref{T:cla:cla}. The isometry type of the universal covering manifold \,$\wt{M}$\, of \,$M$\, (and therefore the local isometry type of \,$M$\,)
is easily determined via the theorem of \textsc{Cartan/Ambrose/Hicks} by computing the restriction of the curvature tensor \,$R$\, of \,$Q$\, to  \,$\liem'$\,.

However, we want to know more: namely the exact global structure of \,$M$\, and the position of \,$M$\, in \,$Q$\,. For this we need to construct totally
geodesic isometric embeddings of suitable Riemannian manifolds into \,$M$\, explicitly, at least for one example per type of Lie triple system.
We will do this below  for all types of
curvature invariant subspaces \,$U$\, except for the type \,$(\Atyp)$\,. In this way we will prove the following table:
\begin{center}
\begin{tabular}{|c|c|c|c|}
\hline
type of \,$\liem'$\, & with ... & isometry class of \,$M$\, & \begin{minipage}{2cm} \begin{center} {\tiny \,$M$\, complex or \\ totally real? \par} \end{center} \end{minipage} \\
\hline
$(\Geo)$ & & \,$\R$\, or \,$\Sph_{L/2\pi}^1$\, & totally real \\
\hline
$(\Geins,k)$ & $2 \leq k \leq m-1$ & $Q^k$ & complex \\
$(\Gzwei,k_1,k_2)$ & $k_1,k_2 \geq 1,\; k_1+k_2 \leq m$ & $(\Sph^{k_1}_{1/\sqrt{2}} \times \Sph^{k_2}_{1/\sqrt{2}})/\{\pm \id\}$ & totally real \\
$(\Gdrei)$ & & $\CP^1 \times \RP^1$ & neither \\
\hline
$(\Peins,k)$ & $1 \leq k \leq m$ & $\Sph^k_{1/\sqrt{2}}$ & totally real \\
$(\Pzwei)$ & & $Q^1$ & complex \\
\hline
$(\Atyp)$ & & $\Sph^2_{\sqrt{10}/2}$ & neither  \\
\hline
$(\Ieins,k)$ & $1 \leq k \leq \tfrac{m}{2}$ & $\CP^{k}$ & complex \\
$(\Izwei,k)$ & $1 \leq k \leq \tfrac{m}{2}$ & $\RP^{k}$ & totally real \\
\hline
\end{tabular}
\end{center}
Here \,$\Sph^k_r \subset \R^{k+1}$\, denotes the \,$k$-sphere of radius \,$r$\,. 
\,$\CP^{k}$\, is equipped with the Fubini-Study metric of constant holomorphic sectional curvature \,$4$\, as usual, and \,$\RP^{k}$\, is equipped
with a Riemannian metric of constant sectional curvature \,$1$\,.

Although it would be possible to describe the totally geodesic embeddings for the various types of Lie triple systems in the general case
via the standard \CQ-structure of \,$\C^{m+2}$\, (see Example~\ref{E:Q:CQ}(b)), for simplicity's sake we here give only one example per type of Lie triple system.
That the embeddings given below indeed map onto totally geodesic submanifolds of \,$Q$\,
is most easily seen via the fact that the connected components of the common fixed point set of a set of isometries is a totally geodesic submanifold 
(see for example \cite{Kobayashi:1972}, Theorem~II.5.1, p.~59). To verify that the images correspond to the stated types of Lie triple systems,
one has to use the description of these types in
Theorem~\ref{T:cla:cla} and the explicit description of the shape operators \,$A_\eta$\, (which constitute the \CQ-structure \,$\A = \A(Q^m,p)$\,
on \,$T_pQ^m$\,) given in Proposition~\ref{P:Q:AR}(a).

\textbf{Types \,$\boldsymbol{(\Geins,k)}$\, and \,$\boldsymbol{(\Pzwei)}$\,.} Let \,$1 \leq k < m$\,. Then
$$ Q^k \to Q^m,\; [z_0,\dotsc,z_{k+1}] \mapsto [z_0,\dotsc,z_{k+1},0,\dotsc,0] $$
is a totally geodesic isometric embedding. Its image corresponds to a Lie triple system of type \,$(\Geins,k)$\, (for \,$k \geq 2$\,) resp.~\,$(\Pzwei)$\, (for \,$k=1$\,).

\textbf{Types \,$\boldsymbol{(\Gzwei,k_1,k_2)}$\,, \,$\boldsymbol{(\Peins,k)}$\, and \,$\boldsymbol{(\Geo)}$\,.} 
Let \,$0\leq k_1,k_2$\, with \,$1 \leq k_1+k_2 \leq m$\,. Then the map
\begin{align*}
\wt{f}_{k_1,k_2}: \Sph^{k_1}_{1/\sqrt{2}} \times \Sph^{k_2}_{1/\sqrt{2}} & \to Q^m, \\
((x_0,\dotsc,x_{k_1}),\; (y_0,\dotsc,y_{k_2})) & \mapsto [x_0, \dotsc, x_{k_1}, i\cdot y_0, \dotsc, i\cdot y_{k_2}, 0, \dotsc, 0]
\end{align*}
is a totally geodesic isometric immersion and a two-fold covering map onto its image with \,$\wt{f}^{-1}(f(x,y)) = \{\pm(x,y)\}$\,. 
It therefore gives rise to a totally geodesic isometric embedding
$$ f_{k_1,k_2}: (\Sph^{k_1}_{1/\sqrt{2}} \times \Sph^{k_2}_{1/\sqrt{2}})/\{\pm \id\} \to Q^m \; . $$
The image of \,$f$\, corresponds to a Lie triple system of type
\,$(\Gzwei, k_1, k_2)$\, (for \,$k_1,k_2 \neq 0$\,) resp.~of type \,$(\Peins,k_1)$\, (for \,$k_2 = 0$\,).

The type \,$(\Gzwei,1,1)$\, warrants special attention: The Lie triple systems of this type are the Cartan subalgebras of \,$\liem$\,, and therefore
the corresponding totally geodesic submanifolds are the maximal flat tori of \,$Q$\,. They can be described in the following way:

We abbreviate \,$r:=\tfrac{1}{\sqrt{2}}$\,  and consider the normal geodesics
\begin{align*}
\wt{\gamma}_1 : \R \to \Sph_r(\R^{m+2}),\; & t \mapsto r\,(\cos(\tfrac{t}{r}) \,,\, 0 \,,\, r\,\sin(\tfrac{t}{r}) \,,\, 0 \,,\, 0 \,,\, \dotsc \,,\, 0) \\
\text{and}\quad
\wt{\gamma}_2 : \R \to \Sph_r(\R^{m+2}),\; & t \mapsto r\,(0\,,\, \cos(\tfrac{t}{r}) \,,\, 0 \,,\, r\,\sin(\tfrac{t}{r}) \,,\, 0 \,,\, \dotsc \,,\, 0)
\end{align*}
and the map
$$ g: \C \to Q,\; t + is \mapsto \pi(\,\wt{\gamma}_1(t) + J\wt{\gamma}_2(s)\,) \; . $$
\,$g$\, is an isometric covering map onto the maximal flat torus \,$\bbT$\, of \,$Q$\,
with \,$z := [1,i,0,\dotsc,0] \in \bbT$\, and \,$T_z\bbT = \pi_*(\R e_3 + i\R e_4)$\, (where \,$e_k$\, is the \,$k$-th canonical basis vector of \,$\C^{m+2}$\,).
The deck transformation group of \,$g$\, is given by the translations in \,$\C$\, by the elements of the lattice
\begin{equation}
\label{eq:tgsub:G2P1:maxtori:Gamma}
\Gamma := \Z\,\tfrac{\pi}{\sqrt{2}}(1+i) \oplus \Z\,\tfrac{\pi}{\sqrt{2}}(1-i) \; .
\end{equation}
%Moreover, we have
%\begin{equation}
%\label{eq:tgsub:G2P1:maxtori:ftang}
%\forall \tau,\sigma \in \R \; : \; T_0f(\tau + i\sigma) = \tau\,v_x + \sigma\,Jv_y \; .
%\end{equation}
It follows that the maximal tori of \,$Q$\, are isometric to the torus \,$\C / \Gamma \cong \Sph^1_{1/2} \times \Sph^1_{1/2}$\,.

Via the maximal tori of \,$Q$\,, in particular the geodesics of \,$Q$\, can be described. The preceding description of the maximal tori
especially gives a means to investigate which geodesics of \,$Q$\, are closed, and what their period is. One obtains the following results
for the maximal geodesic \,$\gamma_v: \R \to Q$\, with \,$\dot{\gamma}_v(0) = v \in TQ$\, in dependence on the characteristic angle \,$\vi(v) \in [0,\tfrac{\pi}{4}]$\,
introduced in Section~\ref{SSe:cla:rk1}: If \,$\tan(\vi(v))$\, is rational, then \,$\gamma_v$\, is periodic. More precisely, the minimal period \,$L$\,
of \,$\gamma_v$\, then is
\begin{enumerate}
\item for \,$\vi(v)=0$\,: \,$L=\sqrt{2}\cdot \pi$\,.
\item for \,$\tan \vi(v) = \tfrac{n_1}{n_2}$\, with \,$n_1,n_2 \in \N$\, relatively prime and \,$n_1,n_2$\, both odd: 
\,$L= \tfrac{\pi}{\sqrt{2}}\cdot \sqrt{n_1^2 + n_2^2}$\,.
\item for \,$\tan \vi(v) = \tfrac{n_1}{n_2}$\, with \,$n_1,n_2 \in \N$\, relatively prime and either \,$n_1$\, or \,$n_2$\, even: 
\,$L= \sqrt{2}\cdot \pi\cdot \sqrt{n_1^2 + n_2^2}$\,.
\end{enumerate}
On the other hand, if \,$\tan(\vi(v))$\, is irrational, then \,$\gamma_v$\, is injective and \,$\overline{\gamma_v(\R)}$\, is a maximal flat torus of \,$Q$\,.

Of course, the totally geodesic submanifolds of \,$Q$\, of type \,$(\Geo)$\, are the traces of the unit speed geodesics of \,$Q$\,.

\textbf{Types \,$\boldsymbol{(\Ieins,k)}$\, and \,$\boldsymbol{(\Izwei,k)}$\,.} Let \,$1 \leq k \leq \tfrac{m}{2}$\,. Then the map
$$ \CP^k \to Q^m,\; [z_0,\dotsc,z_k] \mapsto [z_0,\dotsc,z_k,i\,z_0,\dotsc,i\,z_k,0,\dotsc,0] $$
is a totally geodesic isometric embedding. Its image corresponds to a Lie triple system of type \,$(\Ieins,k)$\,. If it is restricted to a totally geodesic
\,$\RP^k$\, in \,$\CP^k$\,, one obtains another totally geodesic isometric embedding; the image of the latter embedding corresponds to a Lie triple
system of type \,$(\Izwei,k)$\,. 

Note that the totally geodesic submanifolds of \,$Q$\, of type \,$(\Ieins,k)$\, are in fact \,$k$-dimensional complex projective subspaces of the
ambient projective space \,$\CP^{m+1}$\,, and therefore also totally geodesic submanifolds of \,$\CP^{m+1}$\,.

\textbf{Type \,$\boldsymbol{(\Gdrei)}$\,.} The Segre embedding \,$\psi: \CP^1 \times \CP^1 \to \CP^3,\;([z_0,z_1],[w_0,w_1]) \mapsto
[z_0w_0, z_0w_1,z_1w_0,z_1w_1]$\, (see for example \cite{Griffiths/Harris:1978}, p.~192)
is an isometric embedding, whose image \,$Q' := \psi(\CP^1\times \CP^1)$\, is holomorphically congruent to the standard complex quadric \,$Q^2$\,
in \,$\CP^3$\,. If we let \,$g: \CP^3 \to \CP^3$\, be the holomorphic isometry with \,$f(Q') = Q^2$\,, and let \,$C \subset \CP^1$\, be the trace
of a closed geodesic in \,$\CP^1$\, (then \,$C$\, is isometric to \,$\RP^1$\,), the map
$$ \CP^1 \times C \hookrightarrow \CP^1 \times \CP^1 \overset{\psi}{\longrightarrow} Q' \overset{g}{\longrightarrow} Q^2 \longrightarrow Q^m \;, $$
(where the last arrow represents the standard embedding for the type \,$(\Geins,2)$\, described above)
is a totally geodesic isometric embedding. Its image corresponds to a Lie triple system of type \,$(\Gdrei)$\,. 

\textbf{Type \,$\boldsymbol{(\Atyp)}$\,.}
In this case it is not so easy to give a totally geodesic embedding. The obvious approach of calculating the image of an Lie triple system \,$\liem'$\, 
of type \,$(\Atyp)$\, under the
exponential map of \,$Q$\, leads to a very complicated formula for the embedding of the corresponding totally geodesic submanifold. This formula
does not provide any geometric insight, and for this reason we do not reproduce it here.

Rather, we prove in a different way that the totally geodesic submanifold \,$M$\, of \,$Q$\, corresponding to \,$\liem'$\, is isometric to
\,$\Sph^2_{\sqrt{10}/2}$\,. As described in Theorem~\ref{T:cla:cla}, there exists \,$A \in \A$\, and an orthonormal system \,$(x,y,z)$\, in \,$V(A)$\,
so that with \,$a := \tfrac{1}{\sqrt{5}}(2x+Jy)$\, and \,$b := \tfrac{1}{\sqrt{5}}(y+Jx+\sqrt{3}Jz)$\,, \,$(a,b)$\, is an orthonormal basis of \,$\liem'$\,.
Via Proposition~\ref{P:Q:AR}(b) one calculates
$$ \g{R(a,b)b}{a} = \tfrac{2}{5} \; . $$
Because the curvature tensor of the Riemannian symmetric space \,$M$\, is parallel,
it follows that \,$M$\, is a space of constant curvature \,$\tfrac{2}{5} = \tfrac{1}{r^2}$\, with \,$r := \tfrac{\sqrt{10}}{2}$\,,
and therefore \,$M$\, is locally isometric to the sphere \,$\Sph^2_r$\,.
Hence \,$M$\, is isometric either to the sphere \,$\Sph^2_r$\,, or to the real projective space \,$\RP^2$\,
equipped with a Riemannian metric of constant sectional curvature \,$\tfrac{2}{5}$\,. To decide between these two cases, we calculate the length of
closed geodesics in \,$M$\,: 
Let \,$v \in \Sph(T_pM)$\, be given. Because \,$M$\, is a complete, totally geodesic submanifold of \,$Q$\,,
the maximal geodesic \,$\gamma_v: \R \to Q$\, of \,$Q$\, with \,$\gamma_v(0) = p$\, and \,$\dot{\gamma}_v(0) = v$\,
runs completely in \,$M$\, and also is a geodesic of \,$M$\,.
We have \,$\vi(v) = \arctan(\tfrac{1}{2})$\,, therefore it follows from the preceding result on geodesics in \,$Q$\,
that \,$\gamma_v$\, is periodic and that its minimal period is \,$\sqrt{2}\cdot \pi \cdot \sqrt{1^2 + 2^2} = \sqrt{10}\cdot \pi = 2\pi r$\,. 
This shows that \,$M$\, is isometric to \,$\Sph^2_r$\,.

\begin{Remark}
The types \,$(\Gdrei)$\, and \,$(\Atyp)$\, of totally geodesic submanifolds are the ones which are missing from \cite{Chen/Nagano:totges1-1977}
and \cite{Chen/Nagano:totges2-1978} as was described in the Introduction.
\end{Remark}

\section{Conclusion}
\label{Se:conclusion}

The present classification of the totally geodesic submanifolds of the complex quadric was made possible by a combination of the general relations between
roots resp.~root spaces of a symmetric space and the roots resp.~root spaces of its totally geodesic submanifolds from Section~\ref{Se:roots}
with specific results concerning the
geometry of the complex quadric, especially the explicit description of its root spaces in Proposition~\ref{P:Q:roots}(b). It seems likely to me that
the same approach can be used to obtain a classification of totally geodesic submanifolds for the other two infinite series of rank 2 Riemannian symmetric
spaces of compact type, namely the complex 2-Grassmannians \,$G_2(\C^n)$\, and the quaternionic 2-Grassmannians \,$G_2(\HH^n)$\,. Also one might be able
to obtain at least some results on totally geodesic submanifolds in Riemannian symmetric spaces of higher rank
via a refinement of the results from Section~\ref{Se:roots}.

\end{document}